\documentclass[12pt]{article}
\usepackage[utf8]{inputenc}
\usepackage[english]{babel}
\usepackage{graphicx}
\graphicspath{{figs/}}
\usepackage{wrapfig}
\usepackage{amsmath}
\usepackage{enumitem}
\usepackage{cmap}
\usepackage{soul}
\usepackage{color}
\usepackage{amsfonts}
\usepackage{csquotes}

\usepackage[backend=biber,style=numeric-comp, sorting=none, maxcitenames=2,maxbibnames=10,isbn=false,doi=false,url=false,giveninits=true]{biblatex}
\newtheorem{tm}{Theorem}

\title{Cascades of Lorenz attractors in the Shimizu-Morioka model}

\author{Kazakov A.$^1$, Koryakin V.$^1$, Safonov K.$^{1,2}$, Shilnikov A.$^{3}$\\[4pt]
$^1$ \small National Research University Higher School of Economics,\\
\small 25/12 Bolshaya Pecherskaya Street, 603155, Nizhny Novgorod, Russia\\
$^2$ \small Lobachevsky State University of Nizhny Novgorod,\\
\small pr. Gagarina 23, 603022, Nizhny Novgorod, Russia\\
$^3$ \small Neuroscience Institute and Department of Mathematics and Statistics,\\
\small Georgia State University, 100 Piedmont Ave. SE, Atlanta, Georgia 30303, USA,\\
\small {\tt kazakovdz@yandex.ru, vkoryakin@hse.ru, safonov.klim@yandex.ru, ashilnikov@gsu.edu} \\[4pt]
}
\date{}
\begin{document}

\maketitle
\begin{abstract}

The Lorenz attractor is the first example of a robustly chaotic non-hyperbolic attractor. Each orbit of such an attractor has a positive top Lyapunov exponent, and this property persists under small perturbations despite possible bifurcations of the attractor. In this paper, we study the boundary of the Lorenz attractor existence region in the Shimizu-Morioka model. As in the classical Lorenz system, a part of the boundary is associated with the curve $l_{A=0}$, where the first tangency between some Lyapunov subspaces occurs along orbits of the attractor. However, in the Lorenz system, the curve $l_{A=0}$ forms the exact boundary of the Lorenz attractor existence region. Beyond this curve, the attractor is not robustly chaotic, although it may be indistinguishable from the Lorenz attractor in simple numerical experiments. In the Shimizu-Morioka model, the curve $l_{A=0}$ is divided into two parts. The Lorenz attractor existence region adjoins $l_{A=0}$ along the first part of this curve, as in the Lorenz system. Near the second part, as we show, the region of the existence of the Lorenz attractor is fractal. We describe two infinite cascades of disjoint subregions with the Lorenz attractor. One cascade occurs along the curve $l_{A=0}$, another -- in the transversal direction. We show that along the cascades, the Lorenz attractor undergoes ``doubling bifurcations'', leading to a complication of its topological structure.
\end{abstract}

{\bf Keywords: Lorenz attractor, pseudohyperbolicity, homoclinic orbit, inclination-flip bifurcation, Shimizu-Morioka system.} 

\newpage

\section{Introduction}
In this paper, we study bifurcations of the Lorenz attractor in the Shimizu-Morioka model
\begin{equation}
    \begin{cases}
        \dot x = y,\\
        \dot y = x - \lambda y - xz,\\
        \dot z = -\alpha z + x^2,
    \end{cases}
\label{eq_SMsys}    
\end{equation}
which was proposed in \cite{shimizu1980bifurcation} as an alternative model to the classical Lorenz system for high Rayleigh numbers. This system also serves as a normal form for certain classes of codimension-3 bifurcations of equilibrium states and periodic orbits. In particular, system \eqref{eq_SMsys} describes dynamics near a triply degenerate equilibrium state (with triple zero eigenvalue) in $\mathbb{Z}^2$-symmetric systems of ODEs \cite{vladimirov1993low,shil1993normal}. This system also captures the dynamics in a small neighborhood of a triply degenerate fixed point with multipliers $(-1, -1, 1)$ \cite{shil1993normal, GOST05}.

\begin{figure}[h!]
\center{\includegraphics[width=0.9\linewidth]{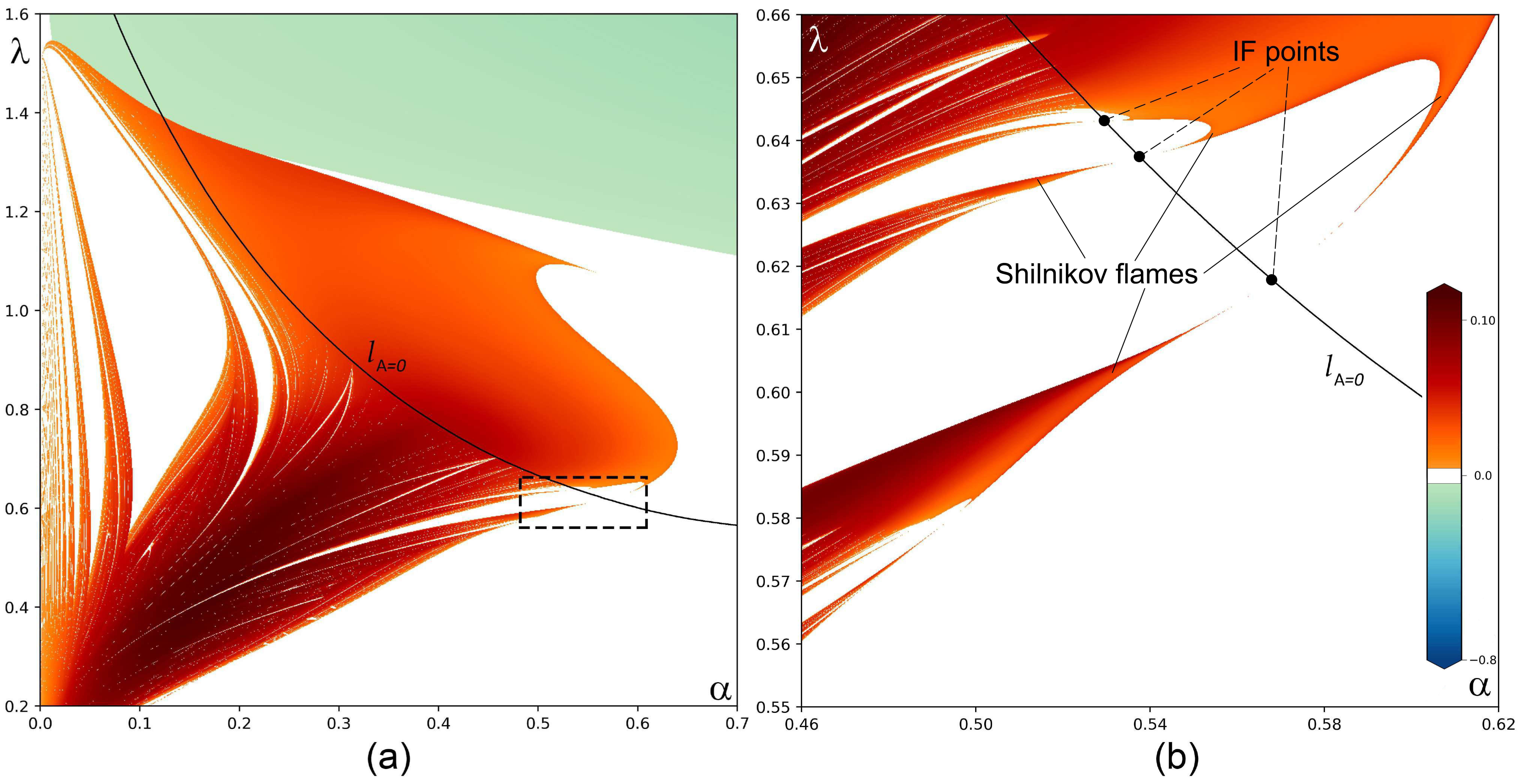}}
\caption{\footnotesize (a) Chart of the top Lyapunov exponent for the Shimizu-Morioka system \eqref{eq_SMsys}. The curve $l_{A=0}$ corresponds to the tangency of the subspaces $E^{ss}$ and $E^{cu}$ along the unstable separatrices of the saddle equilibrium $O(0,0,0)$. The uniform filling of the region above the curve $l_{A=0}$ with orange color indicates the robustness of the attractor. The presence of white stability windows below the curve $l_{A=0}$ indicates that the attractor is not robustly chaotic. (b) An enlarged fragment of the Lyapunov chart shows narrow orange regions (Shilnikov flames) originating from the inclination flip points (IF points) lying on the curve $l_{A=0}$.}
\label{fig_A0}
\end{figure}

An extensive numerical bifurcation analysis of the Shimizu-Morioka system was performed more than 30 years ago in a series of papers \cite{ASh86, ASh89, ASh91, ASh93}. In these papers, the existence of the Lorenz attractors was established for an open region of the parameters, and the main bifurcation curves that bound the existence of the Lorenz attractor were described. In the later paper \cite{barrio2012kneadings}, the authors studied the organization of homoclinic bifurcation curves, which play an important role in understanding the chaotic dynamics associated with Lorenz-like attractors. It was discovered that the region of the existence of the Lorenz attractor includes different disjoint regions -- the so-called \textit{Shilnikov flames} -- originating from codimension-2 points of inclination flip bifurcations (IF points), see Figure~\ref{fig_A0}.

The region of the existence of the Lorenz attractor becomes non-trivial near the curve $l_{A=0}$, where the so-called \textit{separatrix value} vanishes. On this curve, the Lorenz attractor loses its pseudohyperbolicity, which is its principal property. The notion of pseudohyperbolicity was introduced in \cite{TS98, TS08} and characterized by the presence of a continuous splitting of the tangent space in a small neighborhood of the attractor into a direct sum of two invariant linear subspaces $E^{ss}$ and $E^{cu}$. The linearized system restricted to $E^{ss}$ contracts in all directions. In the restriction to $E^{cu}$, the linearized system expands volumes; along some directions, it may not be expanding, but any contraction in this case should be weaker than any contraction in $E^{ss}$, see more details in \cite{GKT21} and Section~\ref{sec_2}. Due to pseudohyperbolicity, the Lorenz attractor possesses the following two important characteristics:
\begin{itemize}
\item each orbit of the attractor has a positive top Lyapunov exponent;
\item this property is robust with respect to small perturbations of the system (e.g., changes of parameters).
\end{itemize}
In other words, pseudohyperbolicity guarantees \textit{robust chaoticity} of an attractor. 

In the Shimizu-Morioka system, the smooth curve $l_{A=0}$ corresponds to the tangency between the subspaces $E^{ss}$ and $E^{cu}$ along the unstable separatrices of the saddle equilibrium $O(0,0,0)$. Most of the region of the existence of the Lorenz attractor lies above the curve $l_{A=0}$, and below this curve, the attractor is mostly not pseudohyperbolic (quasiattractor). In Fig.~\ref{fig_A0}b, we show a part of the parameter plane where the curve $l_{A=0}$ does not properly separate the region of the existence of the Lorenz attractor from stability windows and regions with quasiattractors. For such parameter values, the volume expanding condition may be violated before the tangency between $E^{ss}$ and $E^{cu}$ occurs. This leads to the possible emergence of stable periodic orbits. In particular, we see the white stability windows in the parameter plane intersecting the curve  $l_{A=0}$, see Fig.~\ref{fig_A0}b. The volume expanding condition is preserved only in the \textit{Shilnikov flames} -- narrow regions adjoining the curve $l_{A=0}$. Note that some Shilnikov flames lie below $l_{A=0}$. Thus, the pseudohyperbolic attractor can exist despite the previously occurred tangency between $E^{ss}$ and $E^{cu}$.

In this paper, we study the structure of the region of the existence of the Lorenz attractor near the curve $l_{A=0}$. First, we propose a new numerical method for finding this curve. The method is based on checking pseudohyperbolicity conditions along the unstable separatrices. It helps to locate tangencies between $E^{ss}$ and $E^{cu}$ both inside regions where a chaotic attractor exists and where chaos is transitive. Then, we stress that the curve $l_{A=0}$ consists of two parts. Along the first part, the saddle index $\nu$ of the equilibrium $O$ is less than $1/2$. Recall that the saddle index is the absolute value of the ratio of the leading stable eigenvalue to the unstable one. This part of the curve forms the exact boundary of the existence of the Lorenz attractor. In the second part, the saddle index $\nu(O)$ is in the interval $(1/2,1)$. Near this part, the boundary becomes non-trivial and does not coincide with the curve $l_{A=0}$. We numerically establish that on this part of the curve $l_{A=0}$, there is a cascade of doubling inclination flip bifurcations, and from each inclination flip bifurcation point, two narrow regions of the existence of the Lorenz attractor emerge on both sides of the curve $l_{A=0}$. These regions are bounded by curves of heteroclinic bifurcations where the unstable separatrices of $O$ belong to the stable manifold of saddle periodic orbits.

We show that below the curve $l_{A=0}$, the region of the existence of the Lorenz attractor has a self-similar fractal structure. Namely, besides the curve $l_{A=0}$, we find other similar curves $l^i_{A=0}, i = 2, 4, 8, \dots$ in the parameter plane, which correspond to the tangency between the subspaces $E^{ss}$ and $E^{cu}$ along the unstable separatrices. Then we show that near these curves, the region of the existence of the Lorenz attractor has a similar organization as near the primary curve $l_{A=0}$. We also show that the topological structure of the Lorenz attractor becomes more complex when one moves along the cascade of the tangency curves $l^i_{A=0}, i = 2, 4, 8, \dots$.

The paper is organized as follows. In Section~\ref{sec_2}, we begin with some preliminary notions and results used in this paper. We present historical remarks on the discovery of the Lorenz attractor and give its precise definition, which is based on the pseudohyperbolicity notion and the Afraimovich-Bykov-Shilnikov geometric model. Then, we recall the notion of the kneading invariant, which is a very effective tool in the study of the Lorenz attractor. In this section, we also describe numerical methods which are used for our research. In particular, we explain how to check pseudohyperbolicity, how to detect the curve $l_{A=0}$, and how to build 1D Lorenz maps by a system of ODEs. In Section~\ref{sec_ShilCrit} we discuss the results of L.P. Shilnikov concerning the emergence of the Lorenz attractor via certain codimension-2 homoclinic bifurcations. Such homoclinic bifurcations are key points in the organization of the Lorenz attractor existence region in the Shimizu-Morioka system.

In Section~\ref{sec_scenario}, we describe a one-parameter scenario of the emergence and destruction of the Lorenz attractor. We also propose a new phenomenological scenario leading to the cascade of Lorenz attractors. In Section~\ref{sec_SM}, we study bifurcations of system \eqref{eq_SMsys} in the $(\alpha, \lambda)$-plane and explain the fractal boundary of the region of the existence of the Lorenz attractor. In Section \ref{sec_cascade_SPI}, we show an infinite number of disjoint subregions of the existence of the Lorenz attractor near the curve $l_{A=0}$. In Section~\ref{sec_cascade_SPII}, we stress an infinite number of disjoint subregions of the existence of the Lorenz attractor below the curve $l_{A=0}$. These two sequences of subregions originate from cascades of inclination flip bifurcation points occurring along the curve $l_{A=0}$ and transversely to this curve. We explain this phenomenon and describe the corresponding bifurcation diagrams.

\section{The Lorenz attractor and pseudohyperbolicity} \label{sec_2}

\begin{figure}[h!]
\begin{minipage}[h]{1\linewidth}
\center{\includegraphics[width=1\linewidth]{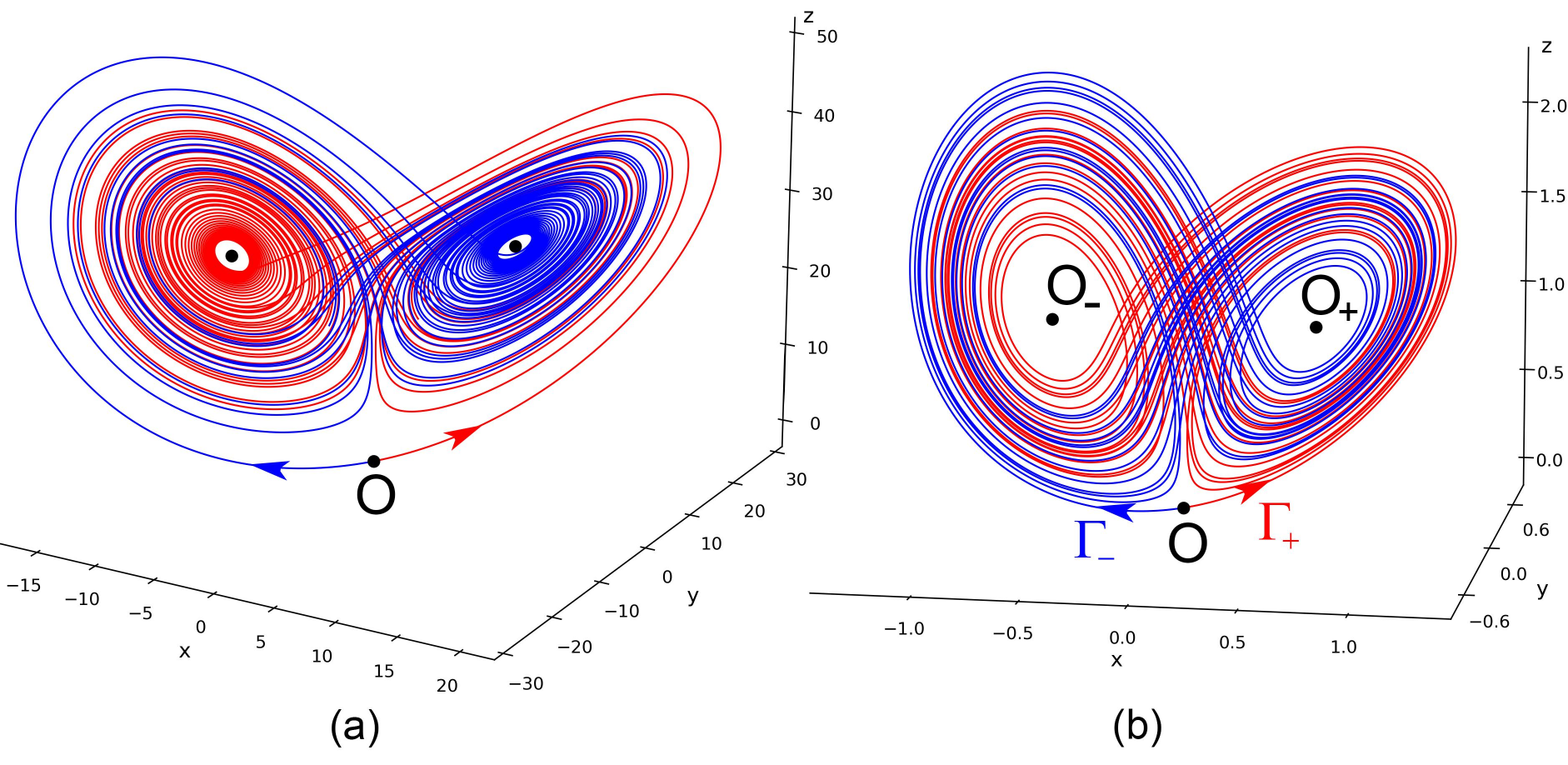}}
\end{minipage}
\caption{\footnotesize Behavior of the unstable separatrices $\Gamma_\pm$ of the saddle equilibrium state $O$ located at the origin: (a) in the Lorenz system \eqref{eq_LorSys} at $b=8/3, \sigma = 10, r=28$; (b) in the Shimizu-Morioka system \eqref{eq_SMsys} at $\alpha=0.4$, $\lambda=0.9$. In both cases, the unstable separatrices draw a chaotic attractor. The points $O^-, O^+$ correspond to saddle-focus equilibrium states.}
\label{fig_LA}
\end{figure}

The Lorenz attractor was discovered by the American mathematician and meteorologist Edward Lorenz in 1963 in the three-dimensional system \cite{Lorenz63}
\begin{equation}
    \begin{cases}
        \dot x = \sigma (y - x),\\
        \dot y = x (r-z) - y,\\
        \dot z = xy - b z.
    \end{cases}
    \label{eq_LorSys}
\end{equation}
This system, as well as system \eqref{eq_SMsys}, is invariant with respect to the axial symmetry
\begin{equation}
\mathcal{S}:\, x \to -x, \ y \to -y, \ z \to z.
\label{eq_LorSym}
\end{equation}
At $b = 8/3$, $\sigma = 10$, and $r = 28$, the Lorenz system demonstrates chaotic dynamics, which can be illustrated by integrating the unstable separatrices of the saddle equilibrium $O(0,0,0)$, see Fig.~\ref{fig_LA}a. Each separatrix rotates around two equilibrium states $O^{\pm}$ located at $\big(\pm \sqrt{b(r-1)}$, $\pm \sqrt{b(r-1)}, r-1\big)$ in an unpredictable stochastic manner. This chaotic behavior can be observed in a large open region of parameters $(b, \sigma, r)$ \cite{shilnikov1980bifurcation, bykov1989boundaries, bykov1992boundaries}.

Similar chaotic dynamics is observed in the Shimizu-Morioka system. For positive $\alpha$, this system also has three equilibrium states: $O(0,0,0)$, which is always a saddle with one-dimensional unstable separatrices $\Gamma_\pm$, and $O_\pm(\pm\sqrt{\alpha},0,1)$, which are either stable or saddle-foci with two-dimensional unstable manifolds. The behavior of the unstable separatrices $\Gamma_\pm$ of the saddle $O$ at $\alpha=0.4$, $\lambda=0.8$ is shown in Fig.~\ref{fig_LA}b. The chaotic attractor in this case is similar to the Lorenz attractor shown in Fig.~\ref{fig_LA}a. Note that for some parameter values, the Lorenz system and the Shimizu-Morioka system are equivalent. Namely, after the following change of coordinates and time (assuming that this change is well-defined)
$$
x_{new}=\frac{\big(1-b/2\sigma\big)^{1/2}}{\big(\sigma(r-1)\big)^{3/4}}\, x,\ \ \ y_{new}=\frac{\big(1-b/2\sigma\big)^{1/2}}{\big(\sigma(r-1)\big)^{5/4}}\, \sigma(y-x), 
$$
$$
z_{new}=\left(z-\frac{x^2}{2\sigma}\right)\, \frac{1}{\sigma(r-1)}, \ \ \ t_{new}=\sqrt{\sigma(r-1)}\, t,
$$
system \eqref{eq_LorSys} takes the form
\begin{equation}
    \begin{cases}
        \dot x = y,\\
        \dot y = x - \lambda y - xz-Bx^3,\\
        \dot z = -\alpha z + x^2,
    \end{cases}
\label{eq_SMsys_ext}    
\end{equation}
which is called the \textit{extended Shimizu-Morioka system}. The parameters of systems \eqref{eq_LorSys} and~\eqref{eq_SMsys_ext} are related by
$$
\alpha=\frac{b}{\sqrt{\sigma(r-1)}}, \ \ \ \lambda=\frac{1+\sigma}{\sqrt{\sigma(r-1)}}, \ \ \ B=\frac{\sqrt{\sigma(r-1)}}{2\sigma-b}.
$$

Nowadays, chaotic attractors of the same shape as the attractors shown in Fig.~\ref{fig_LA} have been found in various nonlinear multidimensional systems, including hydrodynamic models \cite{LZ83, rucklidge1992chaos}, optical systems \cite{pusuluri2018homoclinic, pusuluri2021homoclinic}, normal forms for some classes of local bifurcations \cite{shil1993normal, OT16, gonchenko2022conjoined, karatetskaia2024multi, karatetskaia2024analytic}, etc. 

To study the Lorenz attractor, we need to properly define what we mean by this term. There are several different definitions of the Lorenz attractor \cite{guckenheimer1976strange, guckenheimer1979structural, williams1979structure, ABS77, ABS82, morales2004robust, barros2021up}. In this paper, we follow the approach proposed by Afraimovich-Bykov-Shilnikov \cite{ABS77, ABS82}. In the Afraimovich-Bykov-Shilnikov model, a system has the Lorenz attractor if it generates a singular hyperbolic Poincar\'e map on an appropriate cross-section. We slightly reformulate these assumptions and represent them in a more invariant form using the notion of pseudohyperbolicity \cite{TS98, TS08, CTZ18}.

An $n$-dimensional system $F$ of differential equations is {\it pseudohyperbolic} in a compact forward-invariant domain ${\mathcal D}\subset \mathcal{R}$, if the system possesses the following properties:
\begin{itemize}
\item[(P1)]
at each point $X \in {\cal D}$ there exist two complementary linear subspaces $E^{cu}(X)$ and $E^{ss}(X)$ continuously depending on $X$ such that for all $X \in {\cal D}$ and all $t \geq 0$
$$
D_X F_t E^{cu}(X) = E^{cu}(F_t(X)), \qquad D_XF_t E^{ss}(X) = E^{ss}(F_t(X)),
$$
where $D_XF_t$ is the Jacobian matrix at $X$ of the time-$t$ map $F_t$;
\item[(P2)] the decomposition into $E^{cu}$ and $E^{ss}$ is dominated, i.e., there exist constants $C_1 > 0$ and $\beta > 0$ such that
$$
\|D_XF_t|_{E^{ss}(X)}\| \cdot \|(D_XF_t|_{E^{cu}(X)})^{-1}\| \leq C_1 e^{-\beta t}
$$
for all $t \geq 0$ and $X \in {\cal D}$;
\item[(P3)] the differential $D_XF_t$ restricted to $E^{cu}$ stretches exponentially all $k$-dimensional volumes, where $k={\rm dim}\, E^{cu}$, i.e., there exist constants $C_2 > 0$ and $\sigma > 0$ such that
$$
\det(D_XF_t|_{E^{cu}(X)}) \geq C_2 e^{\sigma t}
$$
for all $t \geq 0$ and $X \in {\cal D}$.
\end{itemize}
The first two conditions ensure that the decomposition into $E^{cu}$ and $E^{ss}$ is robust with respect to small $C^1$-perturbations. The third condition ensures that the top Lyapunov exponent $\Lambda_1$ is positive. Thus, pseudohyperbolicity provides sufficient conditions for chaotic dynamics to be robust with respect to small perturbations.

Below, we consider only the three-dimensional case, in which the definition of the Lorenz attractor is more explicit. Nevertheless, this definition can be naturally generalized to higher dimensions. Assume that a system $F$
\begin{itemize}
\item[(A1)] is invariant with respect to the axial symmetry $\mathcal{S}$;
\item[(A2)] has a symmetric equilibrium $O$ ($\mathcal S(O) = O$) with eigenvalues $\gamma, \lambda_s$, and $\lambda_{ss}$ such that $\gamma>0>\lambda_s>\lambda_{ss}$, and the eigenvector corresponding to $\lambda_s$ coincides with the $z$-axis;
\item[(A3)] is pseudohyperbolic with $\dim E^{cu} = 2$ and $\dim E^{ss} = 1$ in an open forward-invariant domain $\mathcal D$ that contains $O$.
\end{itemize}

\begin{figure}[h!]
\center{\includegraphics[width=0.90\linewidth]{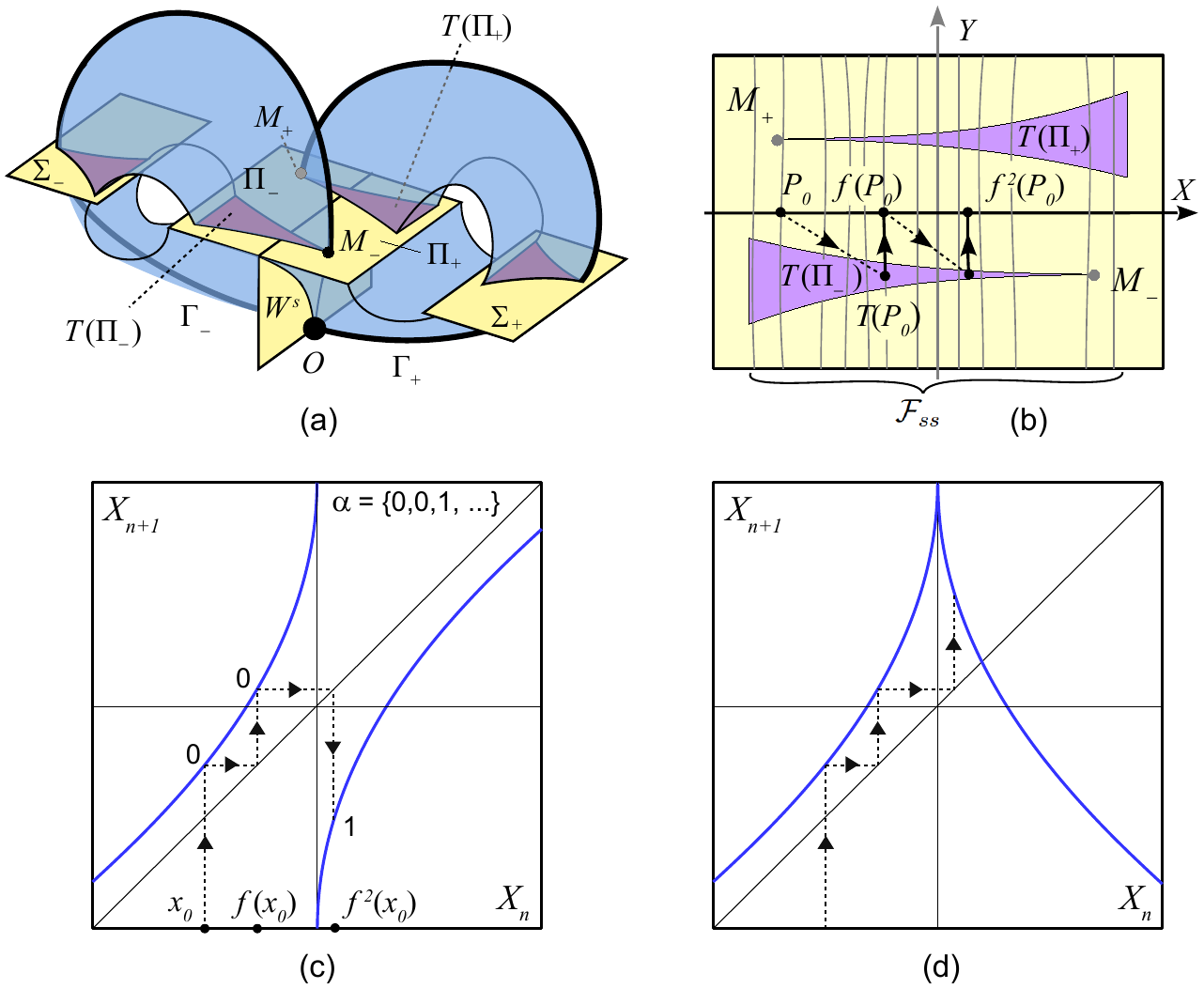}}
\caption{\footnotesize (a) Afraimovich-Bykov-Shilnikov geometric model of the Lorenz attractor. The point $O$ is a saddle equilibrium with a two-dimensional stable invariant manifold $W^s$ and a one-dimensional unstable invariant manifold $W^u$. The manifold $W^u$ is divided by the point $O$ into two orbits -- unstable separatrices $\Gamma_+$ and $\Gamma_-$. The cross-section $\Pi$ is transversal to $W^s$ and is divided by this manifold into two parts $\Pi_+$ and $\Pi_-$. The Poincar\'e map $T$ maps $\Pi_\pm$ to magenta wedges $T(\Pi_\pm)$ with the cusp points $M_\mp$ obtained by the first intersection $\Gamma_\mp$ with $\Pi$. (b) The foliation $\mathcal F_{ss}$ is the projection of the strong-stable invariant foliation provided by the pseudohyperbolicity conditions (P1)--(P3). The iterations $P_0, f(P_0), f^2(P_0),\ldots$ illustrate the action of the factor map $f$ defined by leaves of the foliation $\mathcal F_{ss}$. (c) Graph of the factor map $f$. (d) Graph of the continuous version of the factor map (constructed on the cross-section $\Sigma^+$ or $\Sigma^-$). }
\label{fig_ABS}
\end{figure}

Then, we choose a two-dimensional cross-section $\Pi$ in $\mathcal D$ that transversally intersects the stable manifold $W^s(O)$, see Fig.~\ref{fig_ABS}a. The line $\Pi_0$ of the first intersection $\Pi$ with $W^s(O)$ splits the cross-section into two parts, which we denote by $\Pi_-$ and $\Pi_+$. Further, we assume that:
\begin{itemize}
\item[(A4)] any orbit in $\mathcal D$ either belongs to $W^s(O)$ or intersects $\Pi$;
\item[(A5)] any orbit from $\Pi \backslash \Pi_0$ returns to $\Pi$.
\end{itemize}

Assumptions (A4) and (A5) first describe the butterfly shape of the Lorenz attractor. Second, these assumptions imply that the dynamics of the system $F$ in the domain $\mathcal D$ are completely determined by the Poincar\'e map $T: \Pi_+ \cup \Pi_- \to \Pi$ defined by orbits of the system. The Poincar\'e map is not defined on the line $\Pi_0$, since orbits starting on this line tend to $O$ and do not return to $\Pi$. Nevertheless, we can continuously define the Poincar\'e map from $\Pi_+$ or $\Pi_-$ to the line $\Pi_0$. Namely, let the unstable separatrices $\Gamma_+$ and $\Gamma_-$ intersect the cross-section $\Pi$ at the points $M_+$ and $M_-$, respectively, see Fig.~\ref{fig_ABS}a. Then, the image of a point in $\Pi_+$ and $\Pi_-$ tends to $M_+$ and $M_-$, respectively, when the point tends to $\Pi_0$. Thus, the images of $\Pi_+$ and $\Pi_-$ are wedges with the cusp points $M_+$ and $M_-$. 

Recall that a sequence of points $x_1, \dots, x_N$ is called an $(\varepsilon,\tau)$-orbit if dist$(F_{\tau}(x_j),x_{j+1}) < \varepsilon$ for all $j = 1, \dots, N-1$, where $F_\tau$ is the time-$\tau$ map of the system $F$. A point $y$ is attainable from $x$ if for some $\tau>0$ and any sufficiently small $\varepsilon>0$ there is an $(\varepsilon,\tau)$-orbit connecting the point $x$ with $y$. We define the Lorenz attractor $\mathcal{A}$ as the set of all points attainable from $O$ by $(\varepsilon,\tau)$-orbits. According to~\cite{TS08}, pseudohyperbolicity together with the existence of the cross-section $\Pi$ guarantee the density of $W^s(O)$ in $\mathcal{D}$. This implies that the saddle $O$ is attainable by $(\varepsilon, \tau)$-orbits from any initial point taken in $\mathcal D$. Therefore, the attractor $\mathcal A$ is the only chain-transitive set in $\mathcal D$.


It is important to note that the Lorenz attractor is not hyperbolic (and not structurally stable). For a system of ODEs, hyperbolicity requires the presence of the neutral direction in the tangent decomposition, which is not satisfied for the equilibrium state $O$. Besides, the density of $W^s(O)$ in $\mathcal{D}$ implies the density of the homoclinic butterfly bifurcations in systems with the Lorenz attractor \cite{ABS77, ABS82}. Thus, the attractor changes under small perturbations, but the absorbing domain $\mathcal{D}$ and the whole geometric construction survive. Due to this fact and pseudohyperbolicity, the attractor of the perturbed system remains the Lorenz attractor. According to \cite{TS08}, the chaotic dynamics in systems with the Lorenz attractor also survive small periodic perturbations, which transfer the Lorenz attractor to a \textit{wild pseudohyperbolic attractors} containing homoclinic tangencies \cite{GKT21}.  

\subsection{Numerical verification of pseudohyperbolicity}\label{numerical_verification}

To establish the existence of the Lorenz attractor, we need to verify conditions (A1)--(A5). The Shimizu-Morioka system obviously satisfies symmetry condition (A1). The eigenvalues of $O$ are equal to $-\alpha$ and $-\lambda/2\pm \sqrt{\lambda^2+4}/2$. So, if $0<2\alpha<\lambda+\sqrt{\lambda^2+4}$, then condition (A2) is also fulfilled. Numerical experiments with the Shimizu-Morioka system show that the cross-section $z=1$ satisfies conditions (A4) and (A5) and is suitable for constructing the Poincar\'e map. Thus, to be sure that an observed chaotic attractor is the Lorenz attractor, we need to check the pseudohyperbolicity conditions (A3). The corresponding procedure consists of two steps. We describe them below and illustrate for the attractor shown in Fig.~\ref{fig_LA}b.


The first step is quite simple. It is based on the computation of Lyapunov exponents $\Lambda_1 \geq \Lambda_2 \geq \Lambda_3$ for the attractor. Since points of the attractor are attainable from $O$, we approximate the attractor and Lyapunov exponents by integrating one of the unstable separatrices, e.g., $\Gamma_+$, on a sufficiently long time interval (we take $T = 10^6$). To compute Lyapunov exponents, we apply the standard scheme \cite{BGGS80}. Since $\dim E^{cu} = 2$ and $\dim E^{ss} = 1$ for the Lorenz attractor, the inequality $\Lambda_2 > \Lambda_3$ guarantees the fulfillment of the dominated splitting condition (P2), and the inequality $\Lambda_1 + \Lambda_2 > 0$ ensures the expansion condition (P3). 
For the attractor under consideration we obtain the following set of Lyapunov exponents: $\Lambda_1 = 0.0295, \; \Lambda_2 = 0, \; \Lambda_3 = -1.6095$, i.e., conditions (P2) and (P3) are fulfilled.

The second step is more complicated. It is based on the computation of covariant Lyapunov vectors $V_1(X), V_2(X), V_3(X)$ (these vectors correspond to the Lyapunov exponents $\Lambda_1, \Lambda_2$, and $\Lambda_3$, respectively) and the construction of the subspaces $E^{ss}(X)$ and $E^{cu}(X)$ at each numerically obtained point $X$ of the attractor. To find the covariant Lyapunov vectors, we apply the method proposed by P.~Kuptsov \cite{kuptsov2012fast, KupKuz2018}. The subspace $E^{ss}(X)$ is formed by the vector $V_3(x)$, and the subspace $E^{cu}(X)$ is spanned by the vectors $V_1(x)$ and $V_2(x)$. To establish remaining pseudohyperbolicity condition (P1), we need either to confirm the continuous dependence of $E^{ss}(X)$ and $E^{cu}(X)$ on the point $X$ \cite{GKT21} (straightforward verification of condition (P1)) or to show the absence of tangencies between $E^{ss}(X)$ and $E^{cu}(X)$ \cite{KupKuz2018}. The latter means that the angle $\beta$ between these subspaces is bounded away from $0$ at each point $X$.  According to \cite{GKKT21}, these two approaches are equivalent.

\begin{figure}[h!]
\begin{minipage}[h]{1\linewidth}
\center{\includegraphics[width=1\linewidth]{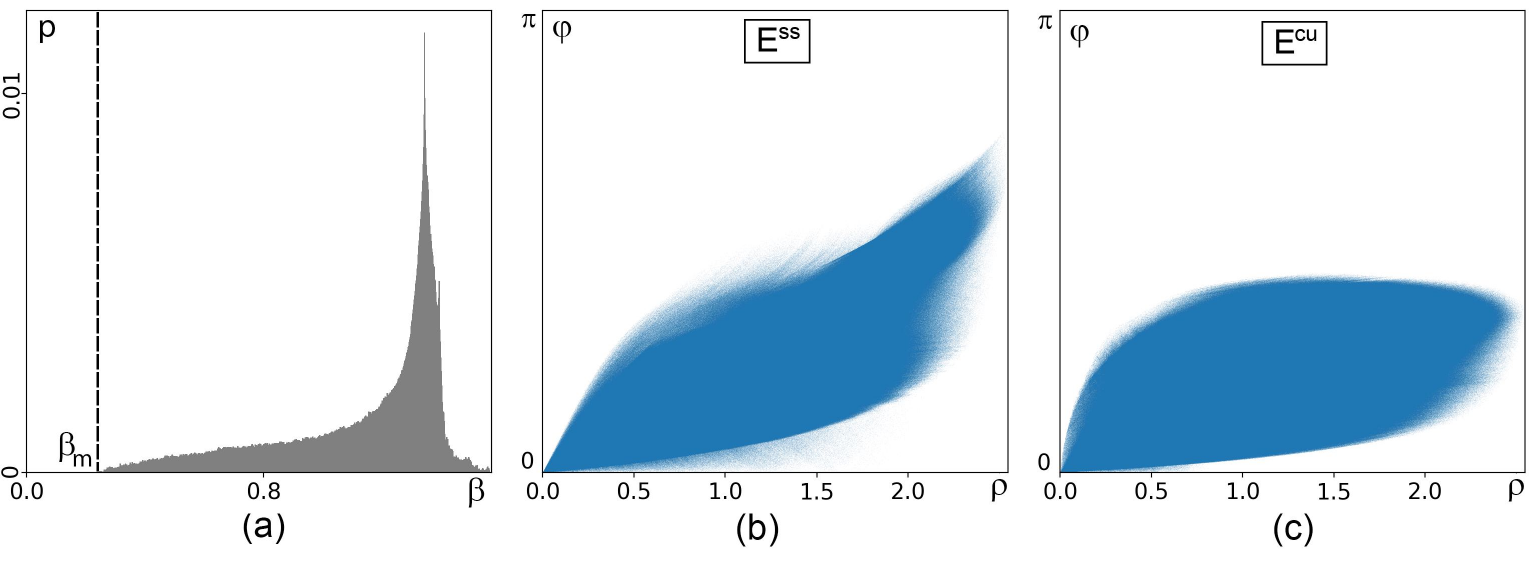}}
\end{minipage}
\caption{\footnotesize Numerical results of the verification of pseudohyperbolicity condition (P1) for the attractor shown in Fig.~\ref{fig_LA}b ($\alpha = 0.4, \lambda = 0.9$). (a) Histogram of the angle between the subspaces $E^{ss}$ and $E^{cu}$. The minimal angle $\beta_m \approx 0.246$ indicates that condition (P1) is satisfied. (b),(c) Continuity diagrams for the subspaces $E^{ss}$ and $E^{cu}$, respectively. In both panels, the cloud of blue points touches the line $\varphi = 0$ only at the point $(0,0)$. This additionally confirms the fulfillment of condition (P1) and also the orientability of the
attractor.}
\label{fig_SM_PhV}
\end{figure}

Figure~\ref{fig_SM_PhV} illustrates the results of the verification of condition (P1) for the attractor shown in Fig.~\ref{fig_LA}b. In panel (a), we plot the histogram of the angles between $E^{ss}(X)$ and $E^{cu}(X)$. To plot it, we calculate the Lyapunov vectors along the unstable separatrix $\Gamma_+$, skipping the initial and final segments of length $10^3$ to ensure that the numerically obtained vectors are aligned along the proper directions. Then, we compute the angle $\beta(X)$ at each point $X$. Here, the minimal angle $\beta_m = 0.246$, i.e., observed chaos is indeed the Lorenz attractor. In addition, in panels (b) and (c), we plot the continuity diagrams for the subspaces $E^{ss}$ and $E^{cu}$, respectively. To plot them, we compute the distance $\rho_{ij}$ between any two points $(X_i,X_j)$ of the attractor and the angle $\varphi_{ij}$ between the corresponding Lyapunov vectors ($V_3(X_i)$ and $V_3(X_j)$ for $E^{ss}$ and $V_1(X_i) \times V_2(X_i)$ and $V_1(X_j) \times V_2(X_j)$ for $E^{cu}$). Condition (P1) is fulfilled if the cloud of points $(\rho_{ij}, \varphi_{ij})$ touches the line $\varphi = 0$ only at the points $(0,0)$ and $(0,\pi)$. If the cloud of points $(\rho_{ij}, \varphi_{ij})$ touches the line $\varphi = 0$ only at the point $(0,0)$, as shown in Figs.~\ref{fig_SM_PhV}b and Figs.~\ref{fig_SM_PhV}c, then the Lorenz attractor is called {\it orientable}. Otherwise, the attractor is called {\it non-orientable}. In Section~\ref{sec_SM}, we show that the non-orientable Lorenz attractor is also observed in the Shimizu-Morioka system. We will also describe bifurcations that change the orientability of the attractor.

The destruction of the Lorenz attractor is associated with the curve $l_{A=0}$ (see Fig.~\ref{fig_A0}), corresponding to the creation of tangency between the subspaces $E^{ss}$ and $E^{cu}$. To find this curve, we use almost the same procedure as for the verification of pseudohyperbolicity condition (P1). The only difference is that we calculate the angles between $E^{ss}$ and $E^{cu}$ along a short piece of the unstable separatrix $\Gamma_+$. Namely, we skip the initial segment of the separatrix of length $0.1$ and consider the next part of the separatrix corresponding to the integration time $t \in [0, 300]$. We use this modification to catch the tangency between $E^{ss}$ and $E^{cu}$ near the first intersection of $\Gamma_+$ with the cross-section $\Pi$.

\begin{figure}[h!]
\begin{minipage}[h]{1\linewidth}
\center{\includegraphics[width=1\linewidth]{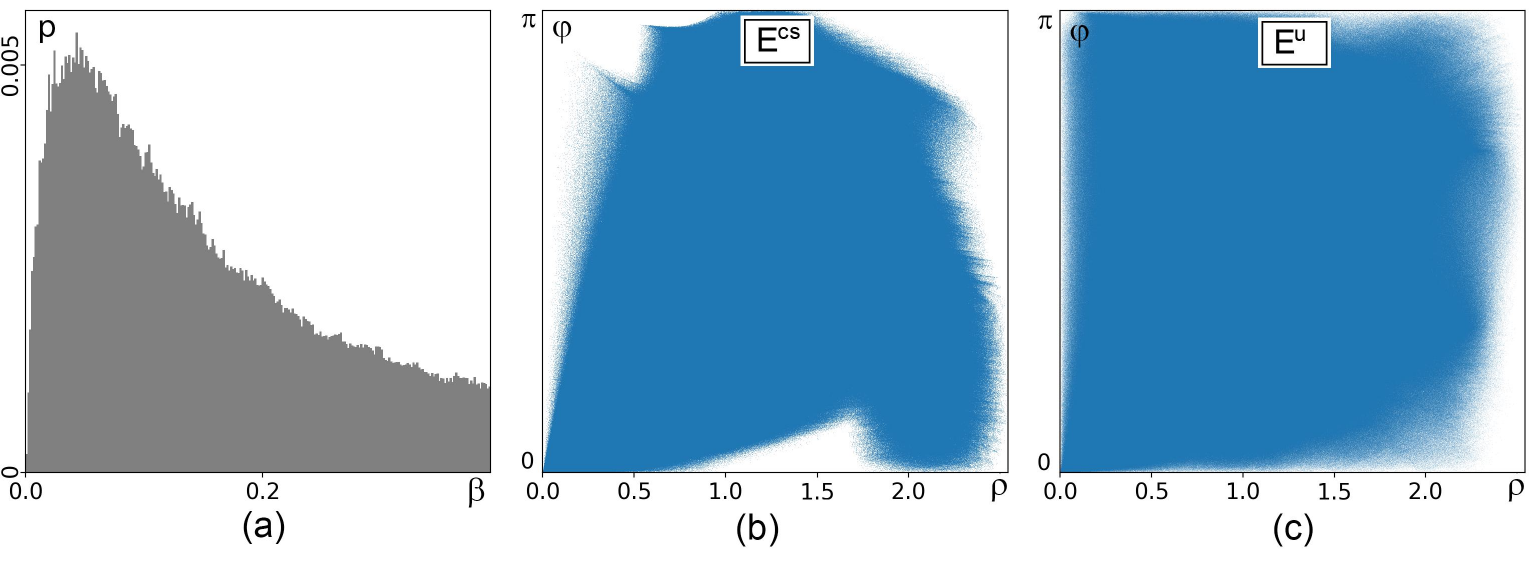}}
\end{minipage}
\caption{\footnotesize Numerical results of the verification of pseudohyperbolicity based on the computation of the subspaces $E^{u}$ (formed by the covariant Lyapunov vector $V_1$) and $E^{cs}$ (spanned by the covariant Lyapunov vectors $V_2$ and $V_3$) for the attractor shown in Fig.~\ref{fig_LA}b ($\alpha = 0.4, \lambda = 0.9$). (a) Histogram of the angle between $E^{u}$ and $E^{cs}$. In this case, the minimal angle between these subspaces vanishes. (b),(c) Continuity diagrams for the subspaces $E^{cs}$ and $E^{u}$, respectively. In panel (c) the cloud of blue points adjoins the line $\varphi=0$, i.e., there is no continuity.}
\label{fig_SM_hypV}
\end{figure}

As noted above, the Lorenz attractor is not hyperbolic. We provide one more argument confirming this fact by calculating the minimum angle between another pair of tangent subspaces: the one-dimensional unstable subspace $E^u(X)$ formed by the vector $V_1(X)$ and the two-dimensional central-stable subspace $E^{cs}(X)$ spanned by the vectors $V_2(X)$ and $V_3(X)$. Figure~\ref{fig_SM_hypV}a shows that this angle vanishes at many points of the attractor. In Figures~\ref{fig_SM_hypV}b and \ref{fig_SM_hypV}c we present the continuity diagrams for the subspaces $E^{cs}$ and $E^u$, respectively. The absence of continuity is clearly visible for $E^u$ where the cloud of points touches the line $\varphi = 0$ at many points.


\subsection{Reduction to the factor map and the kneading invariant} \label{sec_kneadings}

Pseudohyperbolicity conditions imply the existence of the strong-stable invariant foliation in the domain $\mathcal D$. The tangent field of this foliation coincides with the subspace $E^{ss}$. This foliation can be projected onto the strong-stable foliation $\mathcal F_{ss}$ on the cross-section $\Pi$, which is invariant under the Poincar\'e map $T$, see Fig.~\ref{fig_ABS}b. The foliation $\mathcal F_{ss}$ contains the line $\Pi_0$ as a leaf and all its preimages, which densely fill the cross-section. 

The existence of $\mathcal F_{ss}$ allows to factorize the map $T$ over the leaves of this foliation and consider the \textit{factor map} $f: X_n \to X_{n+1}$. The construction of the factor map $f$ is explained in Fig.~\ref{fig_ABS}b. Namely, one takes a point $P_0$ on a smooth curve transversal to the foliation (the line $Y=0$ in Fig.~\ref{fig_ABS}b), and then project its image $T(P_0)$ along the corresponding leaf of $\mathcal F_{ss}$ back onto the curve. The resulting point $f(P_0)$ is the image of $P_0$ by the factor map. The obtained factor map is a piecewise smooth, piecewise monotone, and expanding one-dimensional map with a single discontinuity point, which can be assumed to be fixed at $0$, see Fig.~\ref{fig_ABS}c. Since the foliation is contracting, the dynamics of the factor map completely determine the dynamics of the Poincar\'e map $T$ \cite{malkin1985rotation}. 

Along with the discontinuous factor map $f$, it is convenient to study the continuous map $\tilde{f}$, whose graph is shown in Fig.~\ref{fig_ABS}d. The map $\tilde{f}$ coincides with $f$ for $X<0$ and is equal to $-f$ for $X>0$. The maps $f$ and $\tilde{f}$ are not conjugatd, but there is a natural one-to-one correspondence between orbits of these maps, cf. Fig.~\ref{fig_ABS}c and Fig.~\ref{fig_ABS}d. The continuous map $\tilde{f}$ can also be considered as the factor map of the Poincar\'e map constructed on a different cross-section. Namely, we can take two symmetric cross-sections $\Sigma_-$, $\Sigma_+$ that are intersected by the separatrices $\Gamma_-$, $\Gamma_+$ in the upward direction, and such that  $\Sigma = \Sigma_- \cup \Sigma_+$ satisfies (A4) and (A5), see Fig.~\ref{fig_ABS}a. Then we consider the Poincar\'e map on the cross-section $\Sigma$ and factorize it by the symmetry $S$: if a trajectory intersects $\Sigma_-$, we return to $\Sigma_+$ by the symmetry $S$. The obtained map $\tilde{T}$ maps $\Sigma_+$ to itself, and its factor map is given by a continuous function as shown in Fig~\ref{fig_ABS}d. In this paper, we deal with both Poincar\'e maps $T$ and $\tilde{T}$. The map $\tilde{T}$ is more convenient for numerical experiments, while the map $T$ is more convenient for theoretical studies.

\begin{figure}[h!]
\begin{minipage}[h]{1\linewidth}
\center{\includegraphics[width=1.0\linewidth]{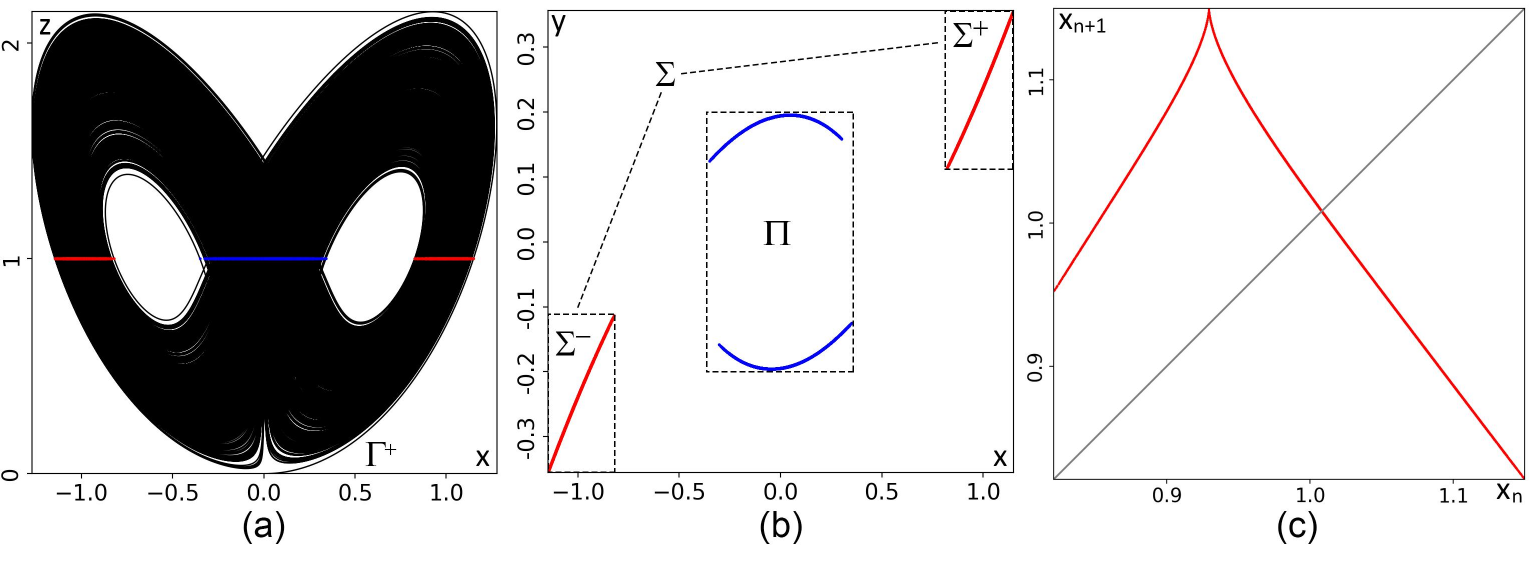}}
\end{minipage}
\caption{\footnotesize Poincar\'e map for the Lorenz attractor shown in Fig.~\ref{fig_LA}b ($\alpha = 0.4, \lambda = 0.9$). (a) A long piece of the unstable separatrix $\Gamma_+$ draws the attractor. (b) The attractor on the secant plane $z=1$. (c) The graph of the 1D Poincar\'e map constructed from the points $(x_n, |x_{n+1}|)$, where $x_n$ is taken in $\Sigma_+$, and $x_{n+1}$ is defined by the Poincar\'e map $(x_{n+1}, y_{n+1}) = T(x_n, y_n)$.}
\label{fig_LorMap}
\end{figure}

In Fig.~\ref{fig_LorMap} we build the 2D Poincar\'e map and the 1D Poincar\'e map (a numerical approximation of the factor map) for the Lorenz attractor shown in Fig.~\ref{fig_LA}b. Here, we take the secant plane $z = 1$ passing through the equilibrium states $O_\pm$ (Fig.~\ref{fig_LorMap}a) and compute the unstable separatrix $\Gamma_+$. The points of its downward intersection are colored in blue, which give the attractor on the cross-section $\Pi$; the points of upward intersection are colored in red, which correspond to the cross-section $\Sigma$, see Figs.~\ref{fig_LorMap}a,b. Then, along this orbit, we construct the 1D Poincar\'e map: for each point $(x_n, y_n)$ on the right part of the cross-section $\Sigma$ we find its image $(x_{n+1}, y_{n+1})$ and plot the point $(x_n, |x_{n+1}|)$ on the graph. The resulting graph of $\tilde{f}$ is shown in Fig.~\ref{fig_LorMap}c. 

The reduction to the factor map allows us to use powerful tools from one-dimensional dynamics in studying the structure and bifurcations of the Lorenz attractor. One such tool is the kneading technique and the notion of the kneading sequence. The kneading sequence is most straightforwardly defined for a piecewise increasing discontinuous factor map $f$, as shown in Fig.~\ref{fig_ABS}c. In this case, we define the map $f$ at the point of discontinuity by continuity from $X\to 0+$. Then, we follow the trajectory of $X=0$ and assign the symbol $0$ if the corresponding iteration lies in $X<0$ and assign the symbol $1$ if the corresponding iteration lies in $X\geq 0$. The obtained infinite sequence is called {\it the kneading sequence}. Note that we can define the factor map by continuity from $X\to 0-$. In this case, we obtain the alternative kneading sequence, where all symbols are replaced with their opposites. According to \cite{malkin1985rotation, li2003smooth}, the kneading sequence is a topological invariant for the map $f$ and the system with the Lorenz attractor. For a decreasing factor map $f$, we use the same construction of the kneading sequence. However, we note that the proper definition of the kneading sequence, which is most commonly used in the one-dimensional dynamics, is slightly different in this case. 

An advantage of the kneading sequence is that it can be calculated without explicit construction of the strong-stable foliation and the factor map. We follow the trajectory of the unstable separatrix (for definiteness, we take the right one, $\Gamma_+$), and calculate maximal and minimal values of the $x$-coordinate. To every maximum we assign the symbol $1$, and to every minimum we assign the symbol $0$. Another advantage of the kneading sequence is that it indicates homoclinic bifurcations: the change in some $n$-th symbol in the kneading sequence means the occurrence of an $n$-round homoclinic loop. Due to this, the kneading diagram helps to display bifurcations of the Lorenz attractor. This approach was used in \cite{barrio2012kneadings} and \cite{xing2014symbolic} to visualize homoclinic bifurcation in the parameter plane of the Lorenz and Shimizu-Morioka systems, respectively. In Section~\ref{sec_SM}, we also apply this approach for the visualization of homoclinic bifurcations near the curve $l_{A=0}$.

\section{Shilnikov criteria}\label{sec_ShilCrit}

In general, a rigorous verification of pseudohyperbolic conditions is quite difficult to do, except for some specific cases, for example, when the system is piecewise linear \cite{belykh2019lorenz, belykh2021sliding}. For the Lorenz system \eqref{eq_LorSys} with the classical parameter values, such a verification was done in \cite{tucker1999lorenz} using computer-assisted proof methods. 


An alternative way to prove the existence of the Lorenz attractor is to find certain degenerate homoclinic bifurcations. In \cite{Shil81}, L.P.~Shilnikov described three homoclinic bifurcations whose unfolding leads to the emergence of the Lorenz attractor. For $\cal S$-symmetric systems, these bifurcations have codimension 2. Thus, by checking the bifurcation conditions at only one point of the parameter plane and for only one trajectory (homoclinic loop), we can establish the existence of the Lorenz attractor in an open region of the parameter plane. Besides, these homoclinic bifurcations are also very essential for describing the entire region of the existence of the Lorenz attractor. Here, we recall two such bifurcations that occur in the Shimizu-Morioka system. We also provide a detailed description of their unfolding in the corresponding two-parameter families. Below, we consider a three-dimensional system $F$ satisfying conditions (A1) and (A2).

\begin{figure}[h!]
\begin{minipage}[h]{1\linewidth}
\center{\includegraphics[width=0.6\linewidth]{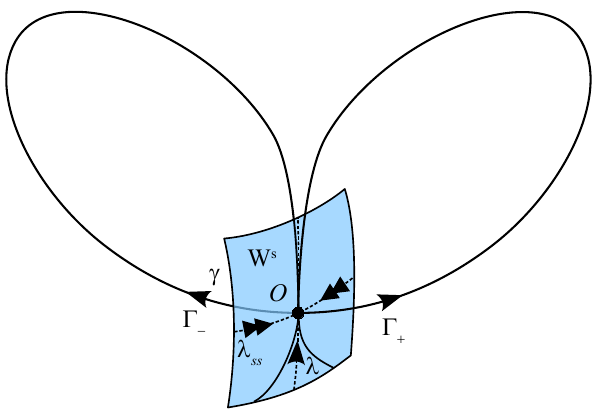}}
\end{minipage}
\caption{\footnotesize Homoclinic butterfly bifurcation: a pair of the unstable separatrices $\Gamma_{\pm}$ return to the saddle $O$ along the leading stable direction.}
\label{fig_hom_but}
\end{figure}

The idea of the Shilnikov criteria comes from the study of the following one-dimensional map
\begin{equation}
X_{n+1} = f(X_n) = \big(-\mu + A|X_n|^{\nu} + \phi(X_n)\big)\, {\rm sign}(X_n).
\label{eq_1DMap}
\end{equation}
Formula \eqref{eq_1DMap} describes the asymptotics of the factor map $f$ near the discontinuity point $X=0$. The parameter $\mu$ is a separatrix splitting parameter. If $\mu=0$, the system has a homoclinic butterfly configuration, i.e., both unstable separatrices $\Gamma_{\pm}$ are double asymptotic to $O$, see Fig.~\ref{fig_hom_but}. The parameter $\nu = -\lambda_s / \gamma$ is the saddle index of the equilibrium. For the existence of the Lorenz attractor, $\nu$ must be less than $1$ (due to assumption A3). The function $\phi(X)$ corresponds to small terms compared to $|Y|^{\nu}$.

\subsubsection*{Shilnikov point I (homoclinic to a neutral saddle)}

Consider a system $F$ that satisfies the following three conditions:
\begin{enumerate}
    \item[(NS1)] the unstable separatrices $\Gamma_{\pm}$ return to $O$ along the stable leading direction ($z$-axis), forming the homoclinic butterfly, as shown in Fig.~\ref{fig_hom_but};
    \item[(NS2)] $\nu = 1$, i.e., the saddle $O$ is neutral;
\end{enumerate}
These two conditions define a codimension-2 bifurcation in the class of $\cal{S}$-symmetric systems. To study this bifurcation, it is natural to embed the system $F$ in a two-parameter family $F_{\mu,\nu}$, depending on the parameters $\mu$ and $\nu$.

For the homoclinic butterfly satisfying assumptions (NS1)--(NS2), we can define the so-called {\it separatrix value}, which describes the dynamics of the system near the homoclinic orbits (the exact definition can be found in \cite{CTZ18,robinson1989homoclinic}). In particular, the separatrix value coincides with the value of the coefficient $A$ in formula \eqref{eq_1DMap} at the bifurcation moment, when $\mu=0$ and $\nu=1$.  We assume that 
\begin{enumerate}
    \item[(NS3)] the separatrix value is non-zero and less than $2$ in absolute value.
\end{enumerate}
In this case, the unfolding of the described bifurcation leads to the emergence of the Lorenz attractor.


\begin{tm}[\cite{Shil81,robinson1989homoclinic,OT16}]
If conditions (NS1)-(NS3) are fulfilled, then in the $(\mu,\nu)$-plane, there exists a region with the Lorenz attractor, containing the bifurcation point $(\mu,\nu) = (0,1)$ on its boundary.
\label{tm_Shil}
\end{tm}

\begin{figure}[h!]
\begin{minipage}[h]{1\linewidth}
\center{\includegraphics[width=1\linewidth]{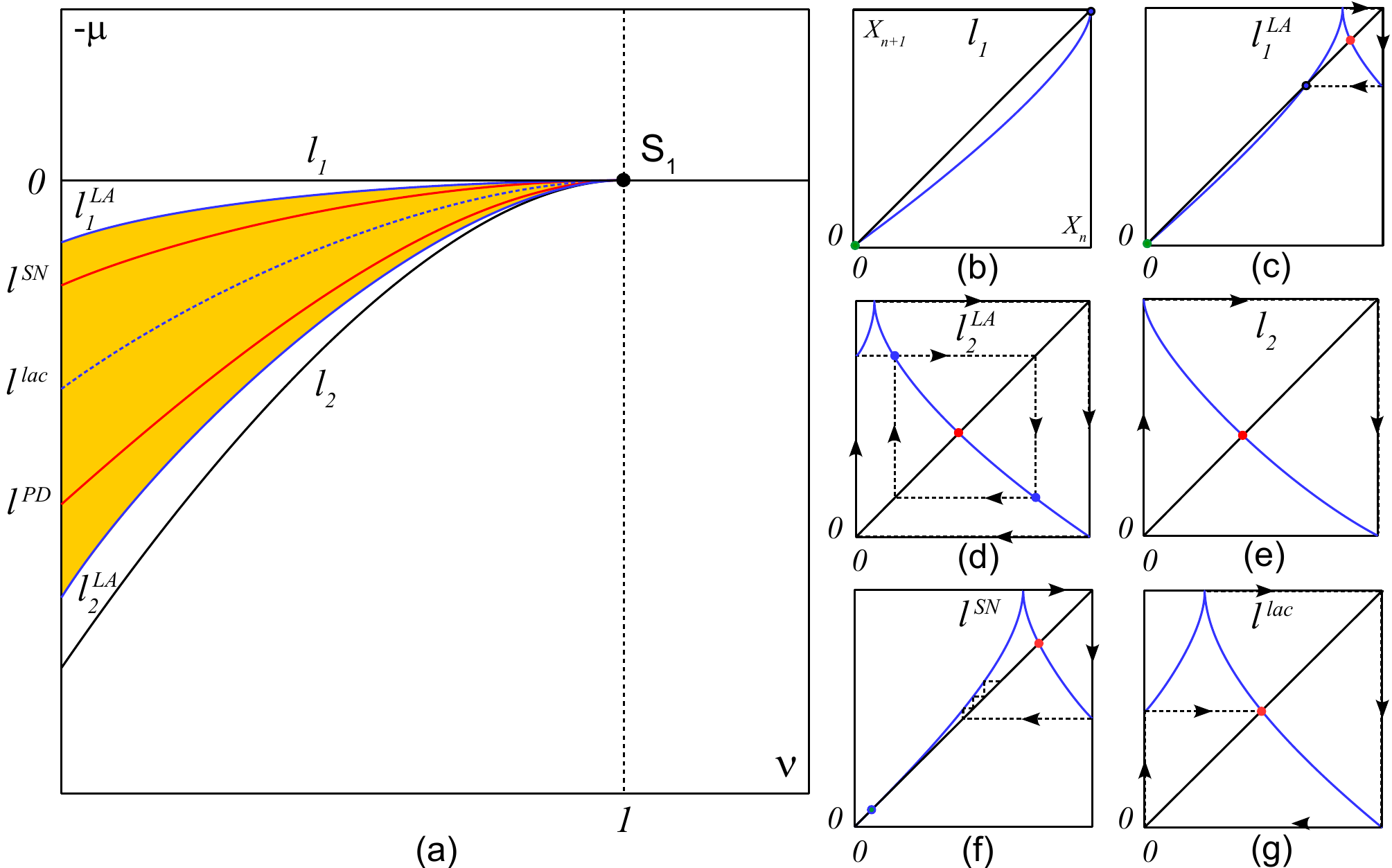}}
\end{minipage}
\caption{\footnotesize (a) Bifurcation diagram for map \eqref{eq_1DMapSigmaTrunc} in the $(\nu,-\mu)$ parameter plane near the Shilnikov point $S_1: (\nu, \mu) = (1,0)$ for $A=0.63$. The region of the existence of the Lorenz attractor (yellow) is bounded by the heteroclinic bifurcation curves $l_1^{LA}$ and $l_2^{LA}$. Panels (b)-(g) illustrate the graphs of map \eqref{eq_1DMapSigmaTrunc} for the bifurcation curves in panel (a). (b) $X=0$ is a fixed point on the line $l_1$ (homoclinic butterfly in the system of ODEs). (c) On the curve $l_1^{LA}$, the second iteration of $X=0$ coincides with an unstable fixed point lying on the increasing branch of the map (heteroclinic connection). (d) On the curve $l_2^{LA}$, the third iteration of $X=0$ coincides with an unstable period-2 point (heteroclinic connection). (e) $X=0$ is a period-2 orbit on the curve $l_2$ (2-round homoclinic butterfly). (f) A tangent (saddle-node) bifurcation occurs on the curve $l^{SN}$. (g) On the curve $l^{lac}$, the third iteration of $X=0$ coincides with the unstable fixed point lying on the decreasing branch of the map (heteroclinic connection). This fixed point does not belong to the attractor below the curve $l^{lac}$, i.e., the attractor has a lacuna. A period-doubling bifurcation occurs on the curve $l^{PD}$.}
\label{fig_ShilCrit_BD}
\end{figure}

To illustrate this theorem, we show in Fig.~\ref{fig_ShilCrit_BD} the bifurcation diagram for the truncated and continuous form of map \eqref{eq_1DMap}:
\begin{equation}
X_{n+1} = \mu - A|X_n|^{\nu}.
\label{eq_1DMapSigmaTrunc}
\end{equation}
In general, the bifurcations diagrams are different for $|A|<1$ and $|A|>1$, see e.g. \cite{kazakov2021bifurcations}. Since for the Shimizu-Morioka system $A \approx 0.63$\footnote{To be more precise, $A \in [0.62597007201516, 0.6267320984754]$, as was rigorously established in \cite{CTZ18}.}, here we consider only the case $0 < A < 1$ (for negative $A$ the diagram is obtained by replacing $\mu\to -\mu$ and $A\to -A$). 

To understand the diagram, let us recall the correspondence between bifurcations of the factor map and the initial system $F$. The orbit starting at $X=0$ describes the behavior of the unstable separatrices of the system $F$. If the point $X=0$ is $n$-periodic, the $n$-round homoclinic loop occurs in the system (in particular, $n=1$ corresponds to the primary 1-round homoclinic loop). If some iteration of the point $X=0$ is an unstable periodic orbit, then a heteroclinic connection occurs: the unstable separatrices lie on the stable manifolds of saddle periodic orbits of the system $F$.

The 1-round homoclinic loop occurs on the line $l_1:\mu = 0$, see Fig.~\ref{fig_ShilCrit_BD}b. The curve $l_1^{LA}: A \mu^{\nu-1} = 2$ corresponds to the heteroclinic connection, when the second iteration of $X=0$ coincides with the unstable fixed point $X=-\mu$, see Fig.~\ref{fig_ShilCrit_BD}c. Near the Shilnikov point $S_1: (\nu, \mu) = (1,0)$, this bifurcation forms the upper boundary of the Lorenz attractor existence region. The bottom boundary here is formed by the curve $l_2^{LA}$, which corresponds to another heteroclinic connection -- when the third iteration of the point $X=0$ coincides with the unstable period-2 orbit $P^2: (p_1, p_2)$, see Fig.~\ref{fig_ShilCrit_BD}d. This period-2 orbit is born from the 2-round homoclinic loop occurring on the curve $l_2: A \mu^{\nu-1} = 1$, see Fig.~\ref{fig_ShilCrit_BD}e. The curves $l_1^{LA}, l_2^{LA}$ and $l_2$ enter the Shilnikov point $S_1$, touching the line $\mu = 0$.

In Fig.~\ref{fig_ShilCrit_BD}a, we also plot several other curves entering the point $S_1$. On the curve $l^{SN}$, the unstable fixed point merges with the stable one due to a tangent (saddle-node) bifurcation, see Fig.~\ref{fig_ShilCrit_BD}f. On the curve $l^{lac}$, the third iteration of the point $X=0$ coincides with the unstable fixed point lying on the decreasing branch of the map, see Fig.~\ref{fig_ShilCrit_BD}g. Below $l^{lac}$, the fixed point becomes unattainable by the orbit of $X=0$ and is no longer part of the Lorenz attractor. In this case, we say that the \textit{attractor has a lacuna}. This unstable fixed point undergoes a period-doubling bifurcation on the curve $l^{PD}$. As a result, the fixed point becomes stable and the unstable period-2 orbit $P^2$ emerges, see Fig.~\ref{fig_ShilCrit_BD}d.\footnote{For map \eqref{eq_1DMap} with $A>0$, a pitchfork bifurcation with a period-2 orbit occurs instead of the period-doubling bifurcation. As a result, a pair of unstable period-2 orbits (one symmetric to the other) emerges.} 

\subsubsection*{Shilnikov point II (inclination flip bifurcation)}

\begin{figure}[h!]
\begin{minipage}[h]{1\linewidth}
\center{\includegraphics[width=0.6\linewidth]{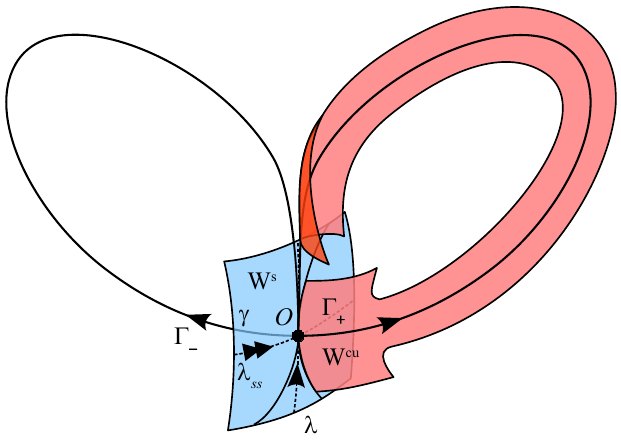}}
\end{minipage}
\caption{\footnotesize Inclination-flip bifurcation: the extended unstable manifold $W^{cu}(O)$ is tangent to the stable manifold $W^s(O)$ along the homoclinic loop.}
\label{fig_inc_flip}
\end{figure}

Now, we consider a system $F$ satisfying:
\begin{enumerate}
    \item[(IF1)] the unstable separatrices $\Gamma_{\pm}$ return to $O$ along the stable leading direction ($z$-axis);
    \item[(IF2)] the extended unstable manifold $W^{cu}(O)$ is tangent to the stable manifold $W^s(O)$ at points of the separatrices $\Gamma_\pm$, see Fig.~\ref{fig_inc_flip};
    \item[(IF3)] $1/2 < \nu < 1$.
\end{enumerate}

As above, the first two conditions defined a codimension-2 bifurcation in the class of $\cal{S}$-symmetric systems. Condition (IF2) is equivalent to the vanishing of the coefficient $A=0$ in map \eqref{eq_1DMap}. Therefore, in this case, it is convenient to consider $\mu$ and $A$ as bifurcation parameters. 

\begin{tm}[\cite{Shil81,rychlik1990lorenz}]
If conditions (IF1)-(IF3) are fulfilled, then in the $(\mu, A)$-plane there exists two disjoint regions with the Lorenz attractor, lying in the half-planes $A>0$ and $A<0$, respectively, and containing the bifurcation point $(\mu, A)=(0,0)$ on its boundary.
\label{tm_Shil2}
\end{tm}

\begin{figure}[h!]
\centering
\includegraphics[width=0.6\linewidth]{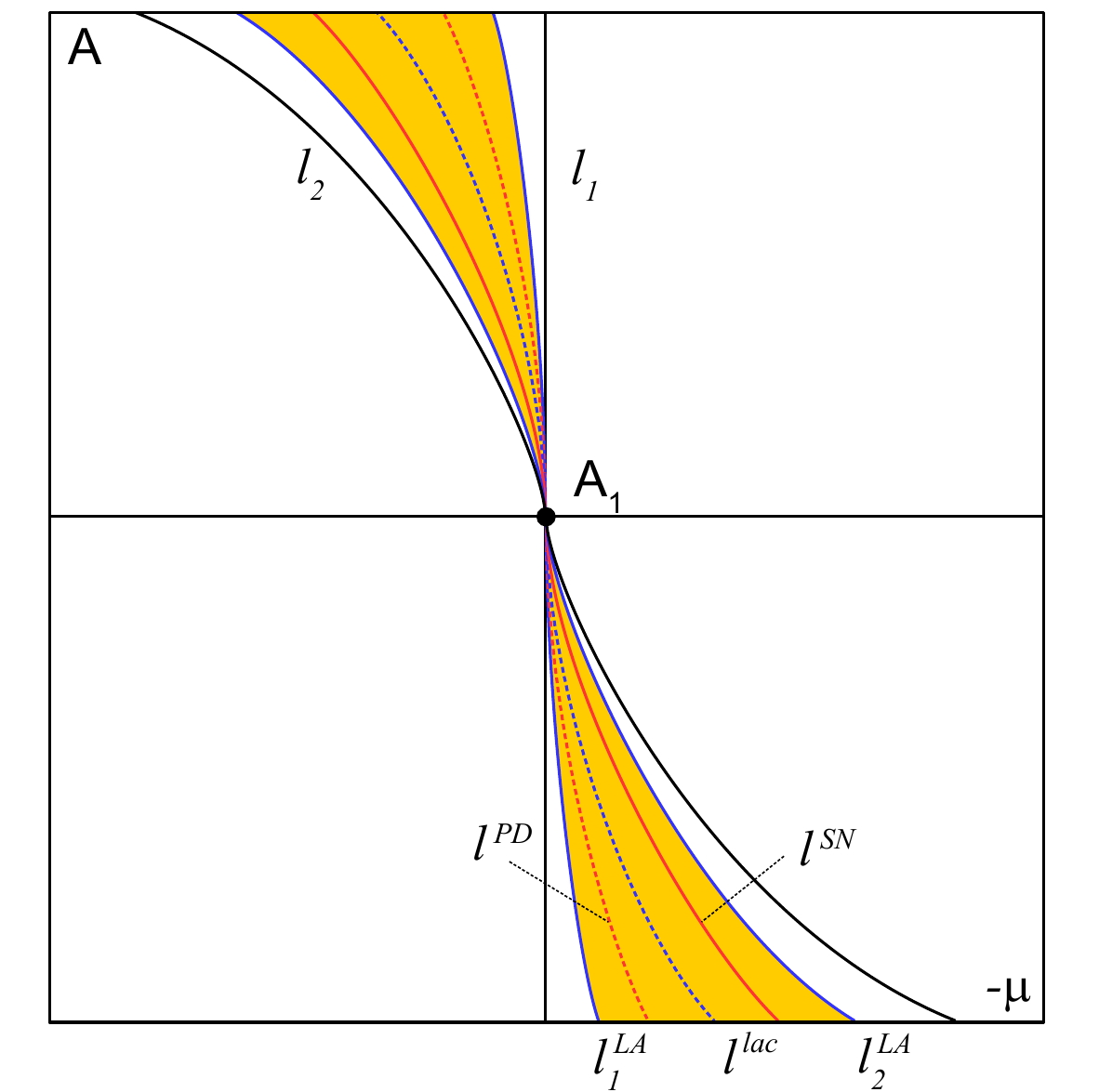}
\caption{\footnotesize Bifurcation diagram for map \eqref{eq_1DMapSigmaTrunc} on the parameter plane $(-\mu, A)$ when $\nu \in (1/2,1)$ near the Shilnikov point $A_1: (\mu,A) = (0,0)$. The meaning of the bifurcation curves is the same as in Fig.~\ref{fig_ShilCrit_BD}.}
\label{fig_ShilCrit_BD2}
\end{figure}

Figure~\ref{fig_ShilCrit_BD2} shows the corresponding bifurcation diagram for map \eqref{eq_1DMapSigmaTrunc}. The Lorenz attractor exists inside the yellow regions bounded by the heteroclinic bifurcation curves $l_1^{LA}$ and $l_2^{LA}$. The meaning of the bifurcation curves is the same as in Fig.~\ref{fig_ShilCrit_BD}.

In the system $F$, changing the sign of $A$ geometrically means changing the topological type of the invariant manifolds $W^{cu}(O)$ and $W^s(O)$. For $A>0$, the primary homoclinic loop (corresponding to the bifurcation curve $l_1$) is orientable: $W^{cu}(O)$ and $W^s(O)$ are annuluses near $\Gamma_{\pm}$. For $A < 0$, this loop becomes non-orientable: $W^{cu}(O)$  and $W^s(O)$ are M\"obius bands near $\Gamma_{\pm}$. On the other hand, the 2-round homoclinic loop (corresponding to the bifurcation curve $l_2$) is always orientable (both for $A < 0$ and $A>0$). As we will show below, this property affects the structure of the Lorenz attractor.

\section{Scenario of the emergence and destruction of the Lorenz attractor} \label{sec_scenario}

In this section, we describe a typical one-parameter scenario of the emergence and destruction of the Lorenz attractor, which, in particular, was observed for both the Lorenz system \cite{shilnikov1980bifurcation} and the Shimizu-Morioka systems \cite{ASh86}. This scenario includes the transition from a system with regular dynamics (a stable equilibrium state or a stable periodic trajectory) to a system with the Lorenz attractor, as well as the further destruction of the attractor's pseudohyperbolicity and the emergence of a Lorenz-like quasiattractor. We also suggest a geometric construction of a one-parameter family of vector fields where the described scenario repeats infinitely many times.


\begin{figure}[h!]
\begin{minipage}[h]{1\linewidth}
\center{\includegraphics[width=1\linewidth]{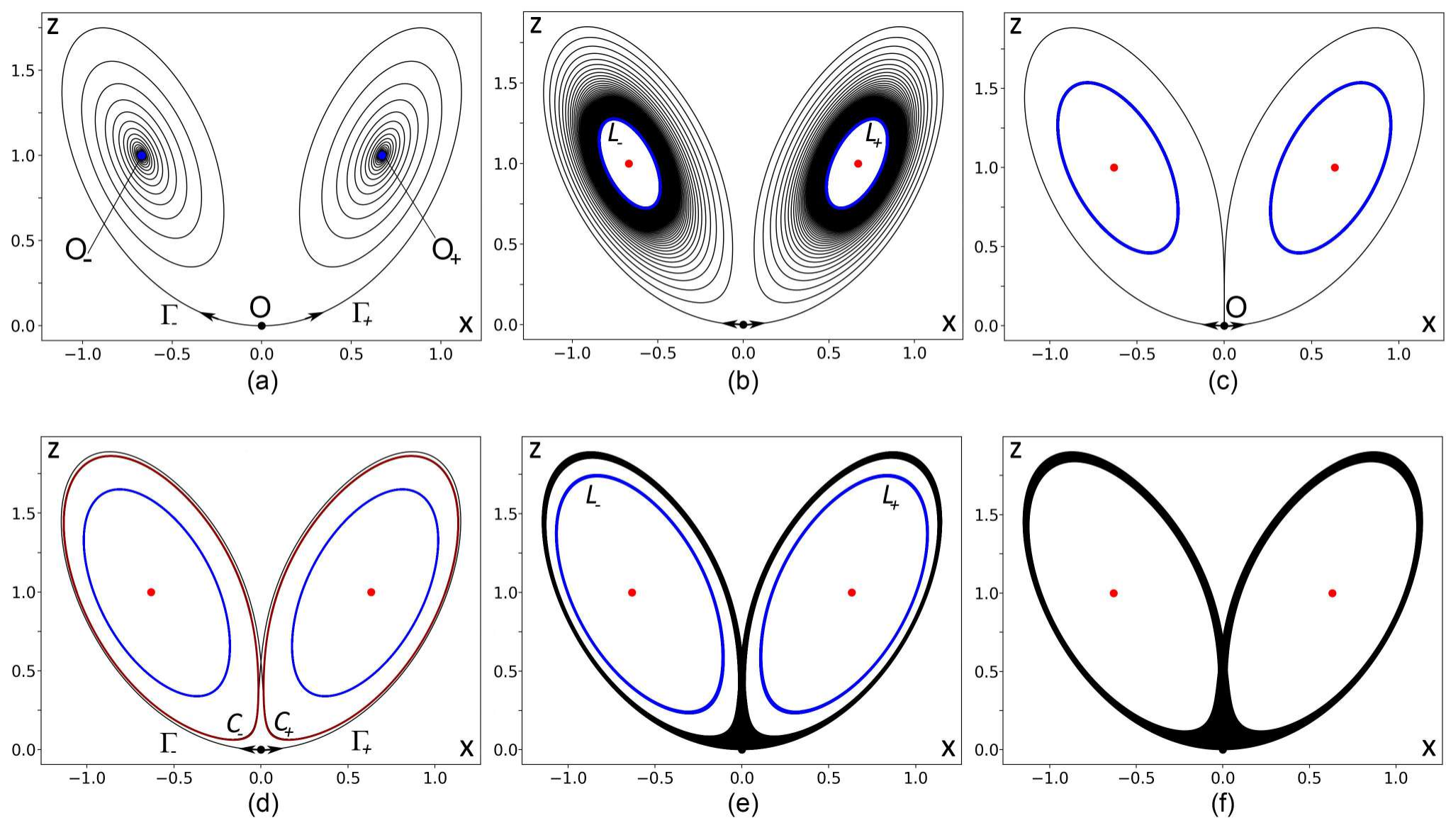}}
\end{minipage}
\caption{\footnotesize Bifurcation scenario of the emergence of the Lorenz attractor in the Shimizu-Morioka system \eqref{eq_SMsys} with fixed $\alpha=0.4$ and $\lambda\in[1.18, 1.35]$. (a) At $\lambda = 1.35$, the unstable separatrices $\Gamma_\pm$ are attracted by the stable equilibria $O_\pm$. (b) At $\lambda = 1.22$, the unstable separatrices $\Gamma_\pm$ tend to the stable periodic orbits $L_\pm$, which are born from the equilibria $O_\pm$ via a supercritical Andronov-Hopf bifurcation. (c) The homoclinic butterfly configuration at $\lambda = 1.205376926$. A pair of periodic saddle orbits $C_\pm$ is born from the homoclinic loops after the splitting of the homoclinic loops. (d) At $\lambda = 1.196740616$, $\Gamma_\pm$ lie on the stable invariant manifolds of $C_\mp$. This bifurcation leads to the emergence of the Lorenz attractor. (e) The phase portrait of the Lorenz attractor at $\lambda = 1.191$. The Lorenz attractor coexists with the pair of stable periodic orbits $L_\pm$. (f) The phase portrait of the Lorenz attractor at $\lambda = 1.18$ after $L_\pm$ merged with $C_\pm$ and disappeared via a saddle-node bifurcation.}
\label{fig_SM_scenario}
\end{figure}

We consider the Shimizu-Morioka system with fixed $\alpha=0.4$ and $\lambda$ running the interval $[0.7, 1.35]$. For $\lambda > \lambda_{AH} \approx 1.234$, the equilibria $O_{\pm}=(\pm\sqrt{\alpha},0,1)$ are stable, and the unstable separatrices $\Gamma_\pm$ of the saddle $O(0,0,0)$ tend to $O_\pm$, respectively, see Fig.~\ref{fig_SM_scenario}a. At $\lambda = \lambda_{AH}$, the equilibria $O_{\pm}$ undergo a supercritical Andronov-Hopf bifurcation. The equilibria $O_{\pm}$ become saddle-foci, and stable periodic orbits $L_\pm$ are born from $O_\pm$, see Fig.~\ref{fig_SM_scenario}b. The unstable separatrices $\Gamma_\pm$ tend to $L_\pm$, respectively.

At $\lambda = \lambda_1 \approx 1.2054$, the unstable separatices $\Gamma_\pm$ return back to $O$ along the $z$-axis, forming the homoclinic butterfly configuration, see Fig.~\ref{fig_SM_scenario}c. When $\lambda$ becomes less than $\lambda_1$, the homoclinic configuration splits and the unstable separatrices rearrange: $\Gamma_+$ tends to $L_-$, while $\Gamma_-$ tends to $L_+$. According to \cite{shilnikov1980bifurcation}, the dynamics of the system become chaotic. The homoclinic bifurcation produces a pair of saddle periodic orbits $C_\pm$ and a nontrivial hyperbolic set $\Omega_0$, the dynamics of which are equivalent to a suspension over the Bernoulli shift with two symbols.

\begin{figure}[h!]
\begin{minipage}[h]{1\linewidth}
\center{\includegraphics[width=1\linewidth]{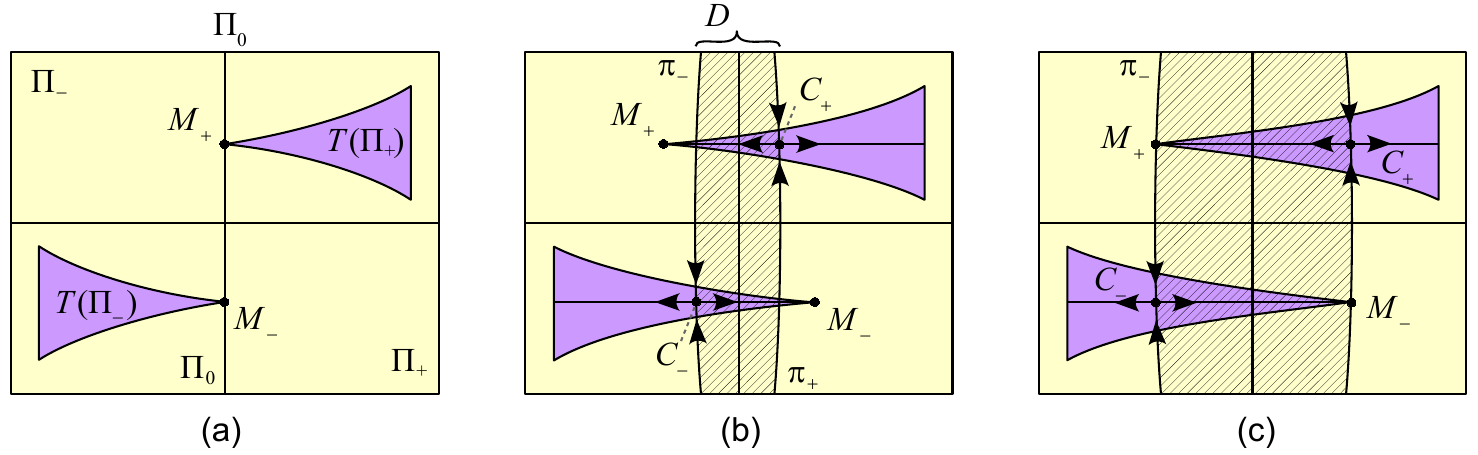}}
\end{minipage}
\caption{\footnotesize The Lorenz attractor emergence in the framework of the Afraimovich-Bykov-Shilnikov geometric model. (a) The homoclinic butterfly bifurcation when $M_{\pm}\in \Pi_0$. (b) Splitting the pair of homoclinic loops creates a pair of periodic saddle orbits $C_\pm$ ($\pi_\pm$ are the traces of their stable manifolds on the cross-section $Pi$). This bifurcation leads to the emergence of a chaotic set, which, in general, is not an attractor, since $dist(\Pi_0, M_\pm) > dist(\Pi_0, \pi_\pm)$. (c) Heteroclinic bifurcation when $M_\pm \in \pi_\pm$.}
\label{fig_GM1a}
\end{figure}

In Figure~\ref{fig_GM1a}a, we show a geometric illustration of the Poincar\'e map at the moment of the homoclinic bifurcation. Recall that the line $\Pi_0$ is the first intersection of the stable manifold $W^s(O)$ and the cross-section $\Pi$. The line $\Pi_0$ splits the cross-section $\Pi$ in two parts $\Pi_{\pm}$. The wedges $T(\Pi_{\pm})$ are the images of $\Pi_{\pm}$ by the Poincar\'e map, and the cusp points $M_{\pm}$ of the wedges are the first intersections of the unstable separatrices $\Gamma_{\pm}$ with $\Pi$. At the moment of homoclinic bifurcation, the points $M_{\pm}$ lie on the line $\Pi_0$, and the wedges $T(\Pi_{\pm})$ do not intersect $\Pi_0$. Note that there are two possibilities for the locations of $T(\Pi_{\pm})$. The first case is when $T(\Pi_+)$ and $T(\Pi_-)$ lie entirely in $\Pi_+$ and $\Pi_-$, respectively, as shown in Fig.~\ref{fig_GM1a}a. This case corresponds to a pair of orientable homoclinic loops. The other case is when $T(\Pi_-)$ lies in $\Pi_+$ and $T(\Pi_+)$ lies in $\Pi_-$. This case corresponds to a pair of non-orientable homoclinic loops. In the Shimizu-Morioka system, the homoclinic loops are orientable at $\alpha=0.4$ and $\lambda\approx 1.2054$. Therefore, we focus our consideration only on the orientable case.

As we mentioned, the homoclinic butterfly splitting leads to the emergence of a pair of periodic saddle orbits $C_-$ and $C_+$. These orbits correspond to saddle fixed points of the Poincar\'e map, see Fig.~\ref{fig_GM1a}b. The lines $\pi_-$ and $\pi_+$ are the traces of the stable manifolds of $C_-$ and $C_+$ on the cross-section $\Pi$. After the homoclinic butterfly splitting, the points $M^{\pm}$ lie outside the domain $\mathcal D$ bounded by the lines $\pi_{\pm}$ (the dashed region in Fig.~\ref{fig_GM1a}). Since the projection of the chaotic set $\Omega_0$ to $\Pi$ consists of orbits that do not leave the domain $\mathcal D$ in the forward and backward iterations of the Poincar\'e map, the chaotic attractor does not exist here.


At $\lambda = \lambda_{LA} \approx 1.1967$, a heteroclinic connection occurs: the unstable separatrices $\Gamma_\pm$ lie on the stable invariant manifolds $W^s(C_\mp)$, respectively, see Fig.~\ref{fig_SM_scenario}d. On the cross-section $\Pi$, the heteroclinic connection corresponds to the situation when $M_{\pm}$ lie on the lines $\pi_{\mp}$, respectively, see Fig.~\ref{fig_GM1a}c. At this moment, forward orbits starting in $\mathcal D$, except for intersections $W^s(O)\cap\Pi$, do not leave this region. An orbit starting in $W^s(O)\cap\Pi$ crosses $\Pi$ a finite number of times, and therefore, the corresponding forward obrit of the Poicnare map is not defined. When $\lambda$ becomes less then $\lambda_{LA}$, the points $M_{\pm}$ lie inside $\mathcal D$, and the region $\mathcal D$ maps strictly inside itself. As a result, the chaotic set $\Omega_0$ becomes attractive (the Lorenz attractor is born) and its absorbing domain is bounded by the stable manifolds of $C_{\pm}$, see Fig.~\ref{fig_SM_scenario}e. Initially, this attractor coexists with the pair of periodic stable orbits $L_\pm$. At $\lambda = \lambda_{SN} \approx 1.1906$, these periodic orbits merge with the saddle periodic orbits $C_\pm$ and disappear via a saddle-node bifurcation, see Fig.~\ref{fig_SM_scenario}f, i.e., the Lorenz attractor $\Omega_0$ becomes the only attractor of the system.

\begin{figure}[h!]
\begin{minipage}[h]{1\linewidth}
\center{\includegraphics[width=1\linewidth]{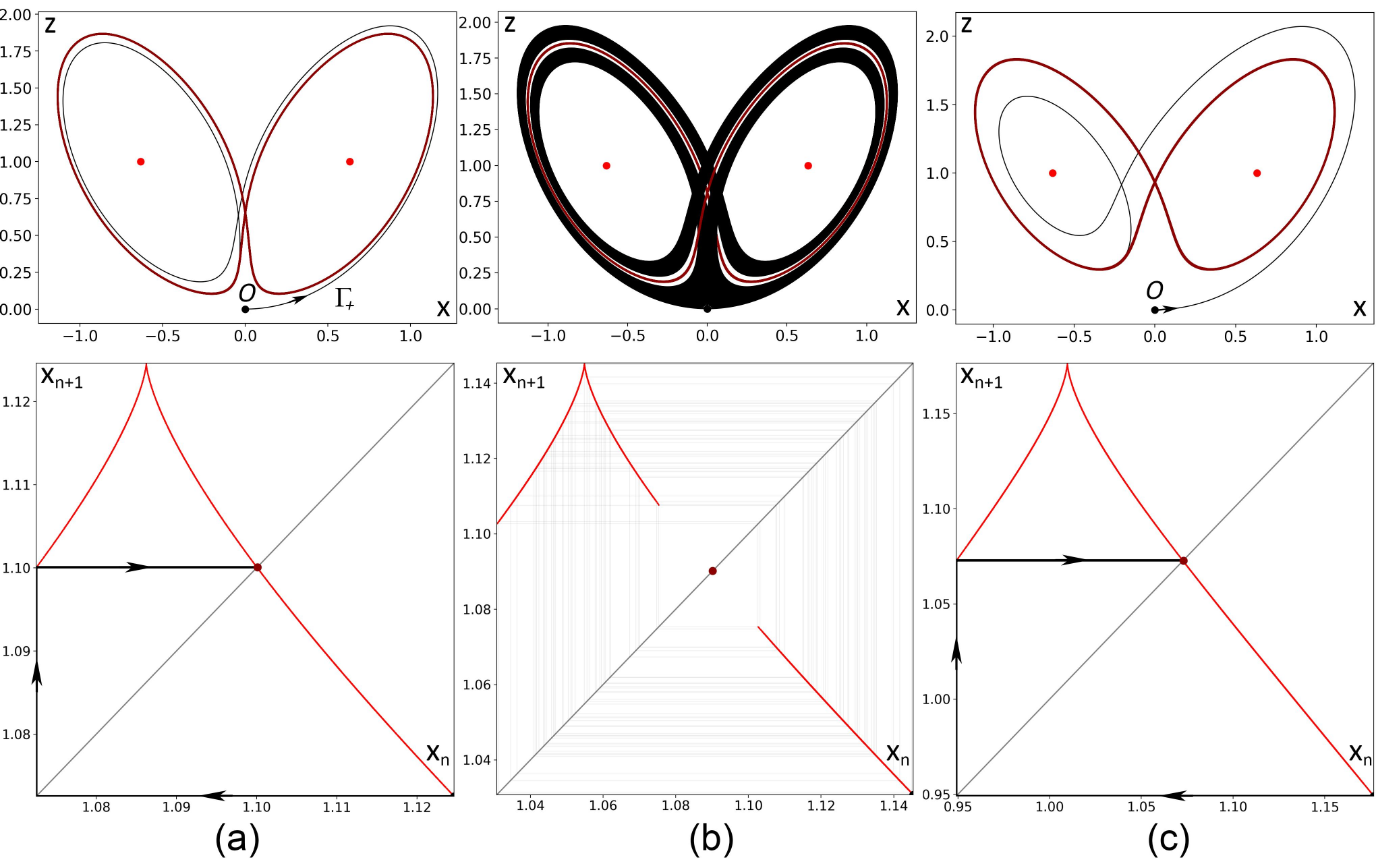}}
\end{minipage}
\caption{\footnotesize The formation and disappearance of a lacuna inside the Lorenz attractor. (a) At $\lambda = 1.1537537538$, the unstable separatrices $\Gamma_\pm$ lies on the stable manifold of the periodic figure-eight saddle orbit. (b) At $\lambda = 1.08$, the periodic saddle orbit is unattainable by orbits of the attractor. The attractor of the 1D map contains two connected components. (c) At $\lambda = 0.9783783784$, $\Gamma_\pm$ again lies on the stable manifold of the periodic saddle orbit; as a result, the lacuna disappears.}
\label{fig_SM_lacuna}
\end{figure}

The next important bifurcation occurs at $\lambda \approx 1.15375$, when the unstable separatrices $\Gamma_\pm$ tend to a periodic figure-eight saddle orbit, see Fig.~\ref{fig_SM_lacuna}a. In the parameter interval $\lambda \in (0.97838, 1.15375)$, this periodic orbit becomes unattainable by orbits of the attractor. The attractor contains a hole, see Fig.~\ref{fig_SM_lacuna}b, which is called a \textit{trivial lacuna} \cite{ABS82}. The lacuna is well observed in the 1D Poincar\'e map, where the attractor becomes two-component. At $\lambda \approx 0.97838$, the backward heteroclinic bifurcation occurs: $\Gamma_\pm$ lie on the same figure-8 periodic orbit and the lacuna disappears. The attractor becomes one-component again, see Fig.~\ref{fig_SM_lacuna}c. The bottom line of Fig.~\ref{fig_SM_lacuna} illustrates these transformations for the numerically constructed 1D Poincar\'e map.


\begin{figure}[h!]
\begin{minipage}[h]{1\linewidth}
\center{\includegraphics[width=1\linewidth]{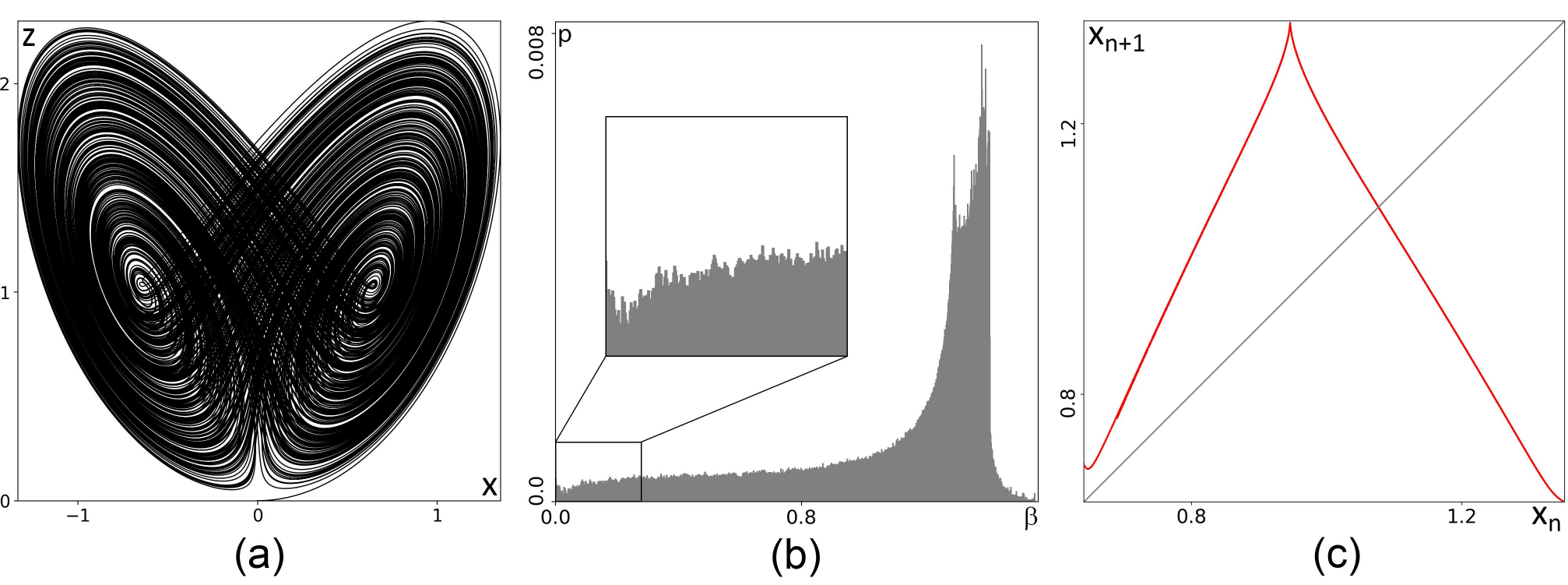}}
\end{minipage}
\caption{\footnotesize Pseudohyperbolicity test of the attractor at $(\alpha, \lambda) = (0.4, 0.76)$. (a) Phase portrait of the attractor. Red points are intersections of the attractor with the cross-section $y=0$. (b) Histogram of the angle between the subspaces $E^{ss}$ and $E^{cu}$. (c) Graph of the numerically constructed 1D map.}
\label{fig_SM_ph_violation}
\end{figure}

At the end of the considered pathway, the Lorenz attractor breaks down due to the violation of pseudohyperbolicity condition (P1). This happens at $\lambda = \lambda_{l_{A=0}} \approx 0.769$. In Figure~\ref{fig_SM_ph_violation} we show the results of the pseudohyperbolicity test for the attractor taken at $(\alpha, \lambda) = (0.4, 0.76)$. Despite the attractor looks like the Lorenz one, see Fig.~\ref{fig_SM_ph_violation}a, it is not pseudohyperbolic due to the tangency between the subspaces $E^{ss}$ and $E^{cu}$ (the minimal angle $\beta_m$ is not bounded away from $0$, see Fig.~\ref{fig_SM_ph_violation}b). In the numerically constructed 1D map, Fig.~\ref{fig_SM_ph_violation}c, we see the emergence of a critical point and a characteristic hook, which also indicate the absence of pseudohyperbolicity. 

\begin{figure}[h!]
\begin{minipage}[h]{1\linewidth}
\center{\includegraphics[width=1\linewidth]{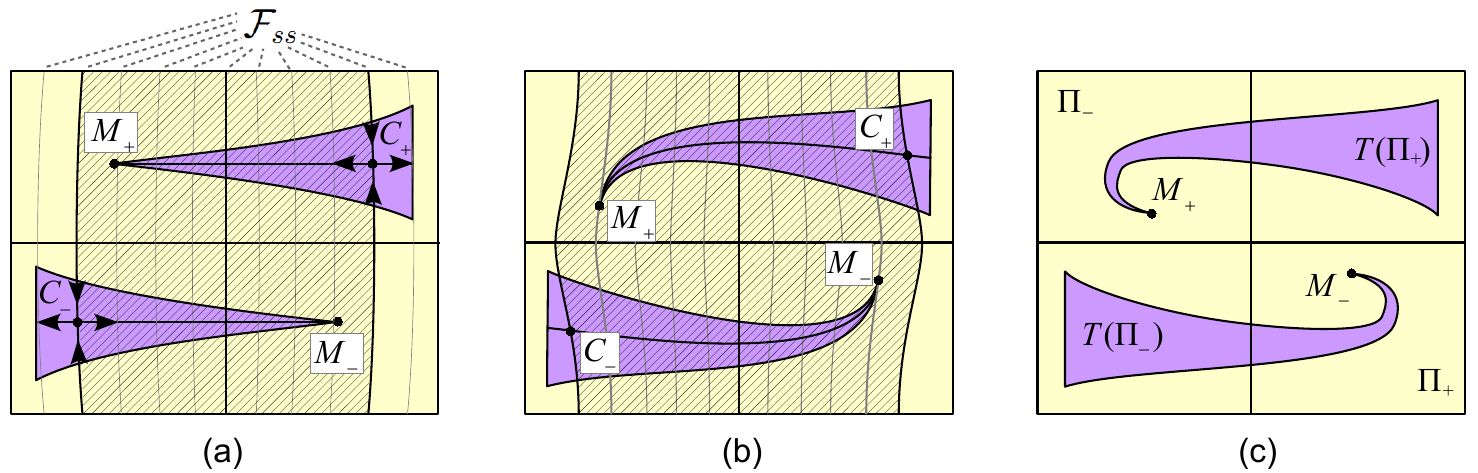}}
\end{minipage}
\caption{\footnotesize The Lorenz attractor destruction in the framework of the Afraimovich-Bykov-Shilnikov geometric model. (a) In the case of the Lorenz attractor, the cross-section $\Pi$ is foliated by the stable invariant foliation $\mathcal F_{ss}$ containing $\pi_\pm$ and $\Pi_0$ as leaves. (b) The moment of destruction of the Lorenz attractor due to the tangency between $\mathcal F_{ss}$ and $T(\Pi_\pm)$. (c) The Poincar\'e map of a Lorenz-like quasiattractor emerging after the pseudohyperbolicity violation.}
\label{fig_GM1b}
\end{figure}

In Fig.~\ref{fig_GM1b} we illustrate the violation of pseudohyperbolicity for the 2D Poincar\'e map. In Fig.~\ref{fig_GM1b}a we show the situation when the attractor is pseudohyperbolic. In Fig.~\ref{fig_GM1b}b the wedges $T(\Pi_\pm)$ become tangent at the points $M_{\pm}$ to some leaves of the foliation $\mathcal F_{ss}$. This corresponds to $A=0$ in the factor map \eqref{eq_1DMap}. The tangency leads to the violation of pseudohyperbolicity conditions (P1): the wedges $T(\Pi_\pm)$ take a hook-shaped form and the foliation $\mathcal F_{ss}$ no longer exists, see Fig.~\ref{fig_GM1b}c.


\subsection{Scenario of an infinite cascade of the Lorenz attractor} \label{sec_cascade}

In the Shimizu-Morioka system, the Lorenz attractor can appear beyond the appearance of the hooks in the Poincar\'e map \cite{ASh89, ASh93, shil1993normal}. A one-parameter scenario of the emergence of such attractor was proposed in \cite{shil1993normal}. The emerging attractor was called \textit{the non-orientable Lorenz attractor}, since its birth is associated with a non-orientable homoclinic loop. In \cite{shil1993normal}, it was also conjectured that such a scenario is observed in the Shimizu-Morioka system. We propose an extended of this scenario -- an infinite cascade of alternating orientable and non-orientable Lorenz attractors together with an infinite cascade of hook-shaped tangencies. In Section~\ref{sec_cascade_SPII}, we give a numerical confirmation that our scenario is observed in the Shimizu-Morioka system.

\begin{figure}[h!]
\begin{minipage}[h]{1\linewidth}
\center{\includegraphics[width=1\linewidth]{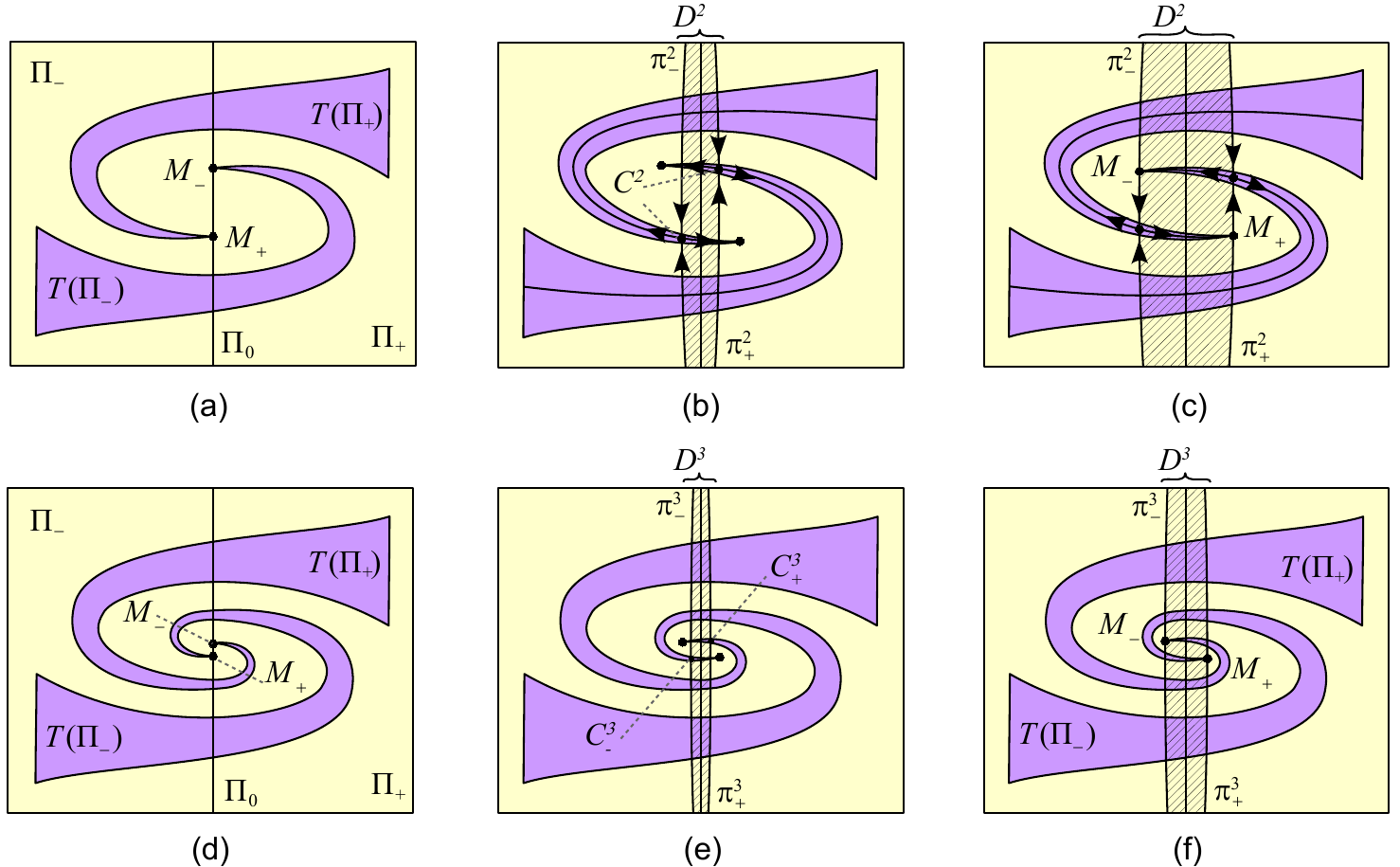}}
\end{minipage}
\caption{\footnotesize Two steps of the scenario of the Lorenz attractor emergence based on the first hooks (top row) and the second hooks (second row) in the framework of the Afraimovich-Bykov-Shilnikov model. (a) Secondary homoclinic butterfly bifurcation. (b) Splitting the homoclinic loops leads to the emergence of a symmetric period-2 saddle orbit $C^2$. The lines $\pi^2_\pm$ are the traces of $W^s(C^2)$. (c) The Lorenz attractor is born when the points $M_\pm$ cross the lines $\pi^2_\pm$ and the region $D^2$ bounded by these curves becomes invariant under the map $T$. (d) The next homoclinic butterfly bifurcation. (e) Splitting the homoclinic loops leads to the emergence of a symmetric pair of saddle periodic orbits $C^3_\pm$. The lines $\pi^3_\pm$ are the traces of $W^s(C^3_\pm)$. (f) The Lorenz attractor is born when the points $M_\pm$ cross the lines $\pi^3_\pm$ and the region $D^3$ bounded by these curves becomes invariant under the map $T$.} 
\label{fig_GM2}
\end{figure}

We begin with the situation shown in Fig.~\ref{fig_GM1b}c, where the tangency between $E^{ss}$ and $E^{cu}$ occurred and the wedges $T(\Pi_{\pm})$ took the hook-shaped form. Let the system changes in such a way that the points $M_{\pm}$ tend again to the line $\Pi_0$, and at some moment, these points lie exactly on $\Pi_0$, see Fig.~\ref{fig_GM2}a. As in Fig.~\ref{fig_GM1a}a, this moment corresponds to the formation of the homoclinic butterfly. However, this homoclinic butterfly is non-orientable, since $T(\Pi_+)$ adjoins the point $M_+$ from the left and $T(\Pi_-)$ adjoins the point $M_-$ from the right. Then, we assume that the system undergoes a similar bifurcation scenario as shown in Figs.~\ref{fig_GM1a} and \ref{fig_GM1b}. 


Namely, the splitting of the non-orientable homoclinic butterfly leads to the emergence of a double-rounded periodic saddle trajectory $C^2$. This periodic trajectory intersects the cross-section $\Pi_0$ at two points $C_{12}, C_{21}$, which forms a symmetric period-2 saddle orbit of the Poincar\'e map $T$, see Fig.~\ref{fig_GM2}b. The lines $\pi^2_-$ and $\pi^2_+$ are the traces of the stable invariant manifold of $C_{12}$ and $C_{21}$, respectively. We consider the region $D^2$ bounded by the lines $\pi^2_-$ and $\pi^2_+$. Right after the splitting of the homoclinic butterfly, $dist(\Pi_0, M_\pm) \gg dist(\Pi_0, \pi^2_\pm)$, i.e., the points $M_{\pm}$ lies outside $D^2$. In particular, the region $D^2$ does not map to itself by the Poincar\'e map $T$. Let the system undergo a heteroclinic bifurcation such that $\Gamma_+$ and $\Gamma_-$ belong to the stable manifold of the double-rounded periodic saddle trajectory $C^2$. Then, the points $M_{\pm}$ lie on the lines $\pi^2_{\mp}$ and the region $D^2$ becomes forward invariant, see Fig.~\ref{fig_GM2}c. As above, the forward invariance relates to all points in $D^2$ except for $\Pi_0$ and its preimages. When the heteroclinic splits and the points $M_{\pm}$ lie inside the region $D^2$, the image of the region $D^2$ lies strictly inside this region. Although pseudohyperbolicity is not fulfilled for all trajectories intersecting $\Pi$, we suppose that pseudohyperbolicity still holds for trajectories intersecting the subregion $D^2$. Then, this region contains a chaotic attractor set, which is the Lorenz attractor.


Pseudohyperbolicity imposes the existence of the stable foliation $\mathcal F^2_{ss}$ on $D^2$ which is invariant by the Poincar\'e map. We assume that the new Lorenz attractor is destroyed via a similar scenario as shown in Fig.~\ref{fig_GM1b}. Namely, the wedges $T(\Pi_\pm)$ become tangent to some leaves of the invariant foliation $\mathcal F^2_{ss}$. After this tangency, the wedges $T(\Pi_\pm)$ get the second hooks, and the foliation no longer exists.


In Figs.~\ref{fig_GM2}d--\ref{fig_GM2}e, we show the next step of our scenario. We suppose that further evolution of the Poincar\'e map leads to another homoclinic butterfly bifurcation, where the points $M_\pm$ again belong to the line $\Pi_0$, see Fig.~\ref{fig_GM2}d. In this case, the homoclinic butterfly is again orientable. Therefore, the splitting of the homoclinic loops leads to the birth of a pair of periodic saddle trajectories $C^3_\pm$, which give a pair of saddle fixed points on $\Pi$. Again, immediately after the splitting, we have $dist(\Pi_0, M_\pm) \gg dist(\Pi_0, \pi^3_\pm)$, where $\pi^3_\pm$ are the traces of $W^s(C^3_\pm)$ on $\Pi_0$. Then, the system undergoes a heteroclinic bifurcation when $M_\pm$ lie on $\pi^3_\pm$, see Fig.~\ref{fig_GM2}f. After this bifurcation, the region $D^3$ bounded by the lines $\pi^3_\pm$ is mapped strictly inside itself by the Poinca\'re map $T$. The region  $D^3$ contains a new Lorenz attractor if the pseudohyperbolicity holds in $D^3$. Then, the attractor is destroyed due to the appearance of new hooks for the wedges $T(\Pi_\pm)$.

By induction, in our scenario, the wedges $T(\Pi_\pm)$ increasingly twist around each other, forming in the limit an infinite spiral. Each even turn of this spiral gives a region with the orientable Lorenz attractor, and each odd turn gives a region with the non-orientable Lorenz attractor. The described phenomenological scenario leads to the infinite sequence of the Lorenz attractors. In Section~\ref{sec_SM}, we show the implementation of this scenario for the Shimizu-Morioka model.



\section{The region Lorenz attractor existence region in the Shimizu-Morioka model} \label{sec_SM}

In this section, we perform a 2-parameter analysis for system \eqref{eq_SMsys} in the $(\alpha,\lambda)$-plane and describe the boundaries of the existence of the Lorenz attractor. The very existence of the Lorenz attractor in this system is a well-known fact \cite{ASh86, ASh89, ASh91, ASh93}. In particular, the fulfillment of the both Shilnikov criteria described in Section~\ref{sec_ShilCrit} was demonstrated in \cite{ASh93}, and the corresponding values of the parameters were found. In \cite{CTZ18}, the fulfillment of the Shilnikov criteria I in system \eqref{eq_SMsys} was established rigorously using computer-assisted proof methods.

\begin{figure}[h!]
\center{\includegraphics[width=1.0\linewidth]{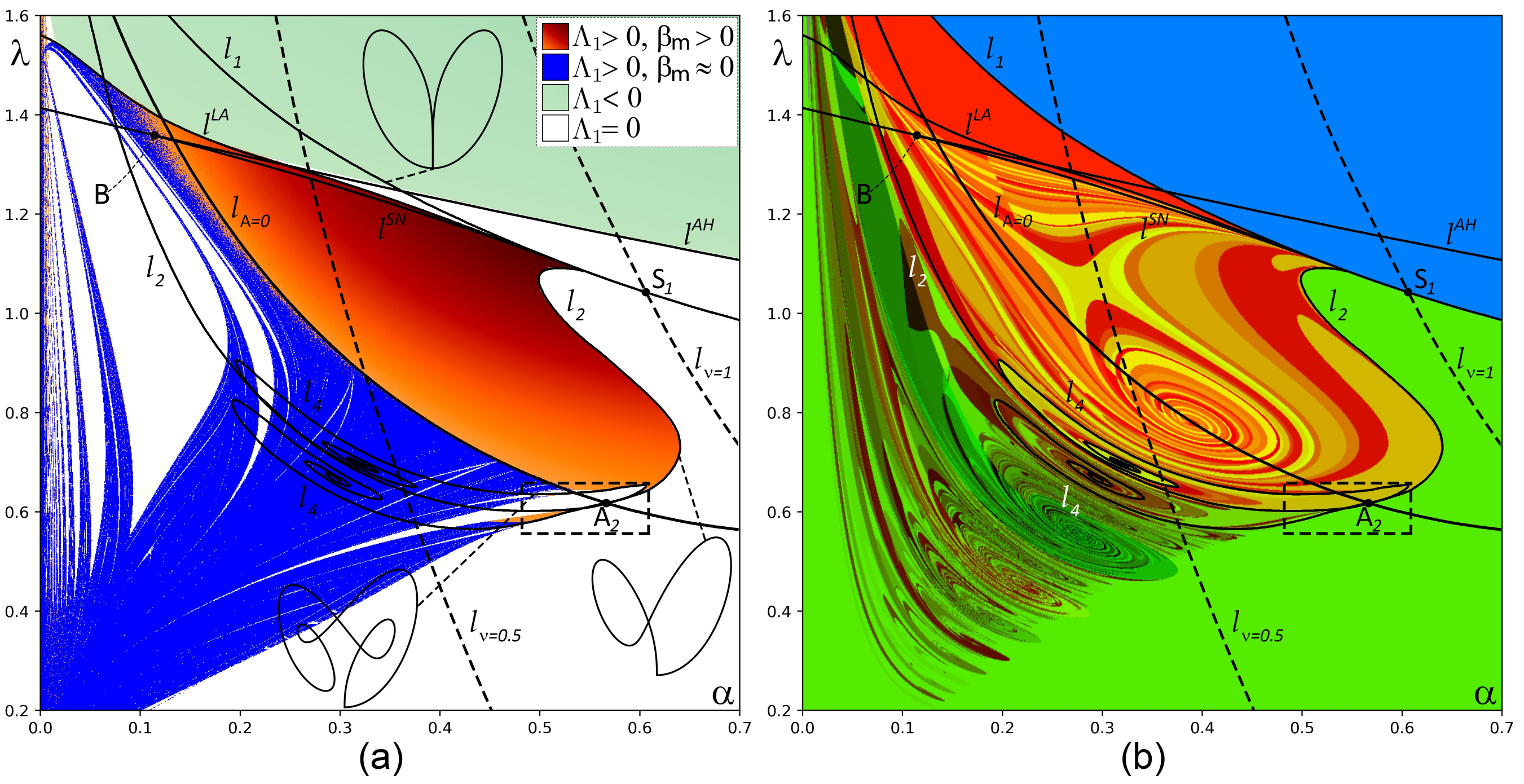}}
\caption{\footnotesize Bifurcation diagram for the Shimizu-Morioka system \eqref{eq_SMsys} superimposed on (a) the chart of the top Lyapunov exponent $\Lambda_1$ and (b) the kneading chart. Pairs of 1-round (primary), 2-round, and 4-round homoclinic loops occur on the curves $l_1$, $l_2$, and $l_4$, respectively. A pair of stable periodic orbits is born as a result of an Andronov-Hopf bifurcation on the curve $l^{AH}$. Then, this pair of stable periodic orbits merges with a pair of saddle periodic orbits (born from the homoclinic loops on $l_1$) as a result of a saddle-node bifurcation on curve $l^{SN}$. The separatrices $\Gamma_\pm$ lie on the stable manifolds of saddle periodic orbits on the curve $l^{LA}$. The curves $l_{\nu=0.5}$ and $l_{\nu=1}$ correspond to the constant values of the saddle index $\nu(O) = 0.5$ and $\nu(O) = 1$, respectively. The minimal angle $\beta_m$ between the subspace $E^{cu}$ and $E^{ss}$ vanishes on the curve $l_{A=0}$. The points $S_1$ and $A_2$ are Shilnikov points of types I and II, respectively.}
\label{fig_SM_Diagr}
\end{figure}

In Figure~\ref{fig_SM_Diagr}, we plot the bifurcation diagram for the Shimizu-Morioka system. We used the MatCont software package \cite{dhooge2008new, de2012interactive} to build most of bifurcation curves. In panel (a), we superimpose the obtained diagram with the colored chart of the top Lyapunov exponent $\Lambda_1$. The color coding for the Lyapunov diagram is depicted in the top-right corner. In panel (b), the same bifurcation diagram is superimposed with the kneading chart. 

The orange and blue colors in panel (a) correspond to the regions with a chaotic attractor ($\Lambda_1 > 0$). In these regions, we additionally check pseudohyperbolicity condition (P1) using the minimal angle method described in Section~\ref{numerical_verification}. In the blue region, the minimal angle $\beta_m$ is not bounded away from zero for some points ($|\beta_m| < \beta^*$, where $\beta^*$ is experimentally chosen to be $0.005$ for Fig.~\ref{fig_SM_Diagr}). Here, pseudohyperbolicity condition (P1) is violated, and there is a tangency between the pair of subspaces $E^{ss}$ and $E^{cu}$. In the orange region, pseudohyperbolicity condition (P1) is satisfied, and we classify the corresponding attractor as the Lorenz attractor. To compute the kneading chart, we use the method described in Section~\ref{sec_kneadings}.
To obtain an informative and noiseless picture, we take $K_N = 15$ and skip the first symbol, since it is always equal to ``1''. Color changes in the kneading chart indicate homoclinic bifurcations. 

Let us explain the type of dynamical regimes in different regions and the meaning of the curves drawn in the colored charts. For the parameter range considered in Fig.~\ref{fig_SM_Diagr}, system \eqref{eq_SMsys} has three equilibrium states: $O$ and $O_\pm$. The equilibrium $O$ is always a saddle with a one-dimensional unstable manifold. The stability of $O_\pm$ changes when crossing the Andronov-Hopf bifurcation curve $l^{AH}$. Above this curve, the equilibria $O_\pm$ are stable and attract the unstable separatrices $\Gamma_\pm$. Below $l^{AH}$, they are of the saddle-focus type. The type of the Andronov-Hopf bifurcation changes at the point $B$ (the Bautin point) corresponding to a degenerate Andronov-Hopf bifurcation. The Andronov-Hopf bifurcation is supercritical to the right of this point and is subcritical to the left of this point. 

Pairs of 1-round, 2-round, and 4-round homoclinic loops occur on the curves $l_1, l_2$, and $l_4$, respectively, see insets in panel (a). The curve $l^{LA}$ corresponds to a heteroclinic connection such that the unstable separatrices $\Gamma_\pm$ lie on the stable manifold of the saddle periodic orbits $C_\pm$, which emerge as a result of the homoclinic bifurcation on the curve $l_1$. On the curve $l^{SN}$, which begins at the point $B$ and ends at the point $S_1$, the saddle periodic orbits $C_\pm$ merge with stable ones and disappear as a result of a saddle-node bifurcation. 

The Lorenz attractor exists below the curve $l^{LA}$ in the orange region, see Fig.~\ref{fig_SM_Diagr}a. Its region of existence originates from the point $S_1$ where the conditions of the Shilnikov criterion I are fulfilled \cite{ASh93, CTZ18}. According to Theorem~\ref{tm_Shil} and Fig.~\ref{fig_ShilCrit_BD}, near the point $S_1$, the upper and lower boundaries of this region are formed by the heteroclinic bifurcation curves $l^{LA}$ and $l^{LA2}$. On the curve $l^{LA2}$, the unstable separatrices $\Gamma_\pm(O)$ lie on the stable manifold of a 2-round saddle periodic orbit, which emerges as a result of the homoclinic bifurcation on the curve $l_2$. Note that in the Shimizu-Morioka system, the curve $l^{LA2}$ passes very close to $l_2$, therefore, we do not plot it in Fig.~\ref{fig_SM_Diagr}.

\subsection{Cascade of inclination flip points along $l_{A=0}$}\label{sec_cascade_SPI}

The lower left boundary of the region of the existence of the Lorenz attractor is associated with the curve $l_{A=0}$, where the first tangency between $E^{ss}$ and $E^{cu}$ occurs. We conjecture that the piece of the curve $l_{A=0}$ lying in the chaotic region ($\Lambda_1 > 0$) and bounded by the inequality $\nu(O) < 1/2$ (to the left of the curve $l_{\nu = 0.5}$) forms the exact boundary of the Lorenz attractor existence region, i.e., this region adjoins the corresponding piece of the curve $l_{A=0}$. 

\begin{figure}[h]
\centering
\includegraphics[width=1\linewidth]{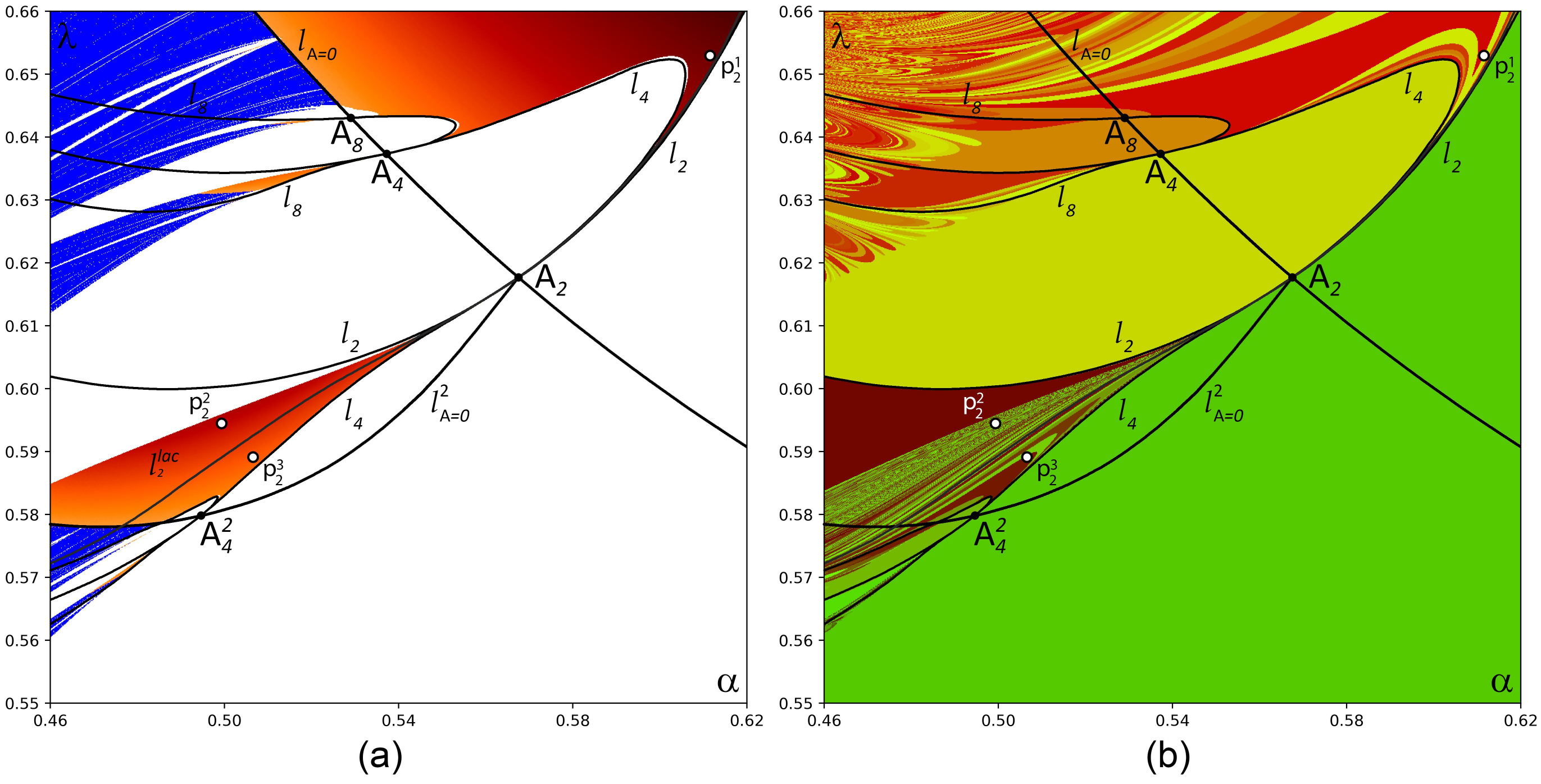}
\caption{\footnotesize Enlarged fragments of (a) the Lyapunov chart and (b) the kneading chart shown in Fig.~\ref{fig_SM_Diagr}. Pairs of 2-round, 4-round, and 8-round homoclinic loops occur on the curves $l_2$, $l_4$, and $l_8$, respectively. These curves intersect the curve $l_{A=0}$ at the points $A_2$, $A_4$, and $A_8$, where the conditions of the Shilnikov criterion II are satisfied. The curve $l_{A=0}^2$ is the curve of the second tangency between $E^{ss}$ and $E^{cu}$. On the curve $l_2^{lac}$, the unstable separatrices $\Gamma_\pm$ lie on the stable manifold of a 4-round symmetric saddle periodic orbit.}
\label{fig_A2}
\end{figure}

Now, we show that the boundary is much more complicated near another piece of the curve $l_{A=0}$ (to the right of the curve $l_{\nu = 0.5}$). In particular, in Fig.~\ref{fig_SM_Diagr}a we see a white region intersecting the curve $l_{A=0}$. This region is located near the point $A_2$, where the curves $l_2$ and $l_{A=0}$ intersect. The point $A_2$ corresponds to the inclination flip bifurcation for the pair of homoclinic loops. At this point, the conditions of the Shilnikov criterion II are fulfilled. Figure~\ref{fig_A2} shows an enlarged fragment of the bifurcation diagram near the point $A_2$. For panel (a), the threshold of the minimal angle $\beta^*$ is experimentally chosen to be equal to $0.0006$. To obtain an informative and noiseless kneading chart in panel (b), we set $K_N = 28$ and skip the first symbol, since it is always equal to ``1''.

The bifurcation diagram near the point $A_2$ is arranged in accordance with Theorem~\ref{tm_Shil2} and the bifurcation diagram shown in Fig.~\ref{fig_ShilCrit_BD2}. In particular, a pair of regions with the Lorenz attractor comes from this point: one above $l_{A=0}$ and the other -- below this curve. Besides, the curve $l_4$ of 4-round homoclinic loops passes through the point $A_2$. Its upper branch intersects $l_{A=0}$ at the inclination flip bifurcation point $A_4$. Similarly, two other regions with the Lorenz attractor and the curve $l_8$ of 8-round homoclinic loops come from this point. The upper branch of the curve $l_8$ intersects $l_{A=0}$ at the next inclination flip bifurcation point $A_8$, and so on. Thus, along the curve $l_{A=0}$, we can observe a cascade of doubling inclination flip bifurcation points. 

\begin{figure}[h!]
\center{\includegraphics[width=0.7\linewidth]{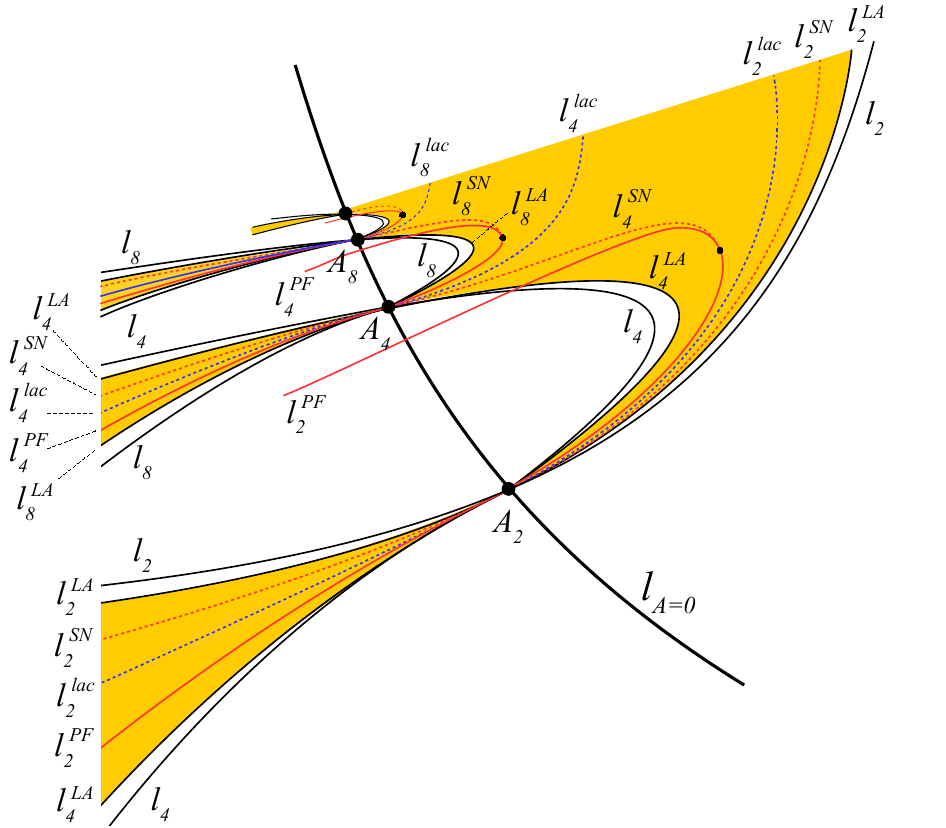}}
\caption{\footnotesize Sketch of the bifurcation diagram near the curve $l_{A=0}$ when $1/2 < \nu < 1$. The Lorenz attractor exists in the yellow region bounded by the heteroclinic bifurcation curves $l_n^{LA}$, where $\Gamma_\pm$ lie on the stable manifolds of $n$-rounded saddle periodic orbits, $n = 2, 4, \dots$. The curves $l_n^{SN}$ correspond to saddle-node bifurcations of periodic orbits. The curves $l_n^{PF}$ correspond to pitchfork bifurcations with symmetric $n$-round periodic orbits; the bifurcations are degenerate in the bold points on these curves. The curves $l_n^{SN}, n > 2$ start from these points and pass through the points $A_n$. The curves $l_n^{lac}$ correspond to the emergence of lacunae on the base of symmetric period-$n$ saddle periodic orbits ($\Gamma_\pm$ lie on the stable manifold of these orbits).}
\label{fig_SM_SketchA0}
\end{figure}

Based on the bifurcation diagram shown in Fig.~\ref{fig_ShilCrit_BD2}, we schematically draw the organization of the region of the existence of the Lorenz attractor along this cascade, see Fig.~\ref{fig_SM_SketchA0}. Above $l_{A=0}$, the region of the existence of the Lorenz attractor has one connected component with the boundary formed by the curves $l_2^{LA}$, $l_4^{LA}, l_8^{LA}, \ldots$. The curves $l_n^{LA}$ correspond to heteroclinic connections when $\Gamma_\pm$ lie on the stable manifold of an $n$-round saddle periodic orbit, $n = 2, 4, \dots$. Below $l_{A=0}$, the region of the existence of the Lorenz attractor has many connected components coming from the points $A_2, A_4, A_8, \ldots$. Each connected component is bounded by a pair of curves $l_{2n}^{LA}$ and $l_{n}^{LA}$, $n = 2, 4, \dots$.

\begin{figure}[h!]
\centering
\includegraphics[width=1.0\linewidth]{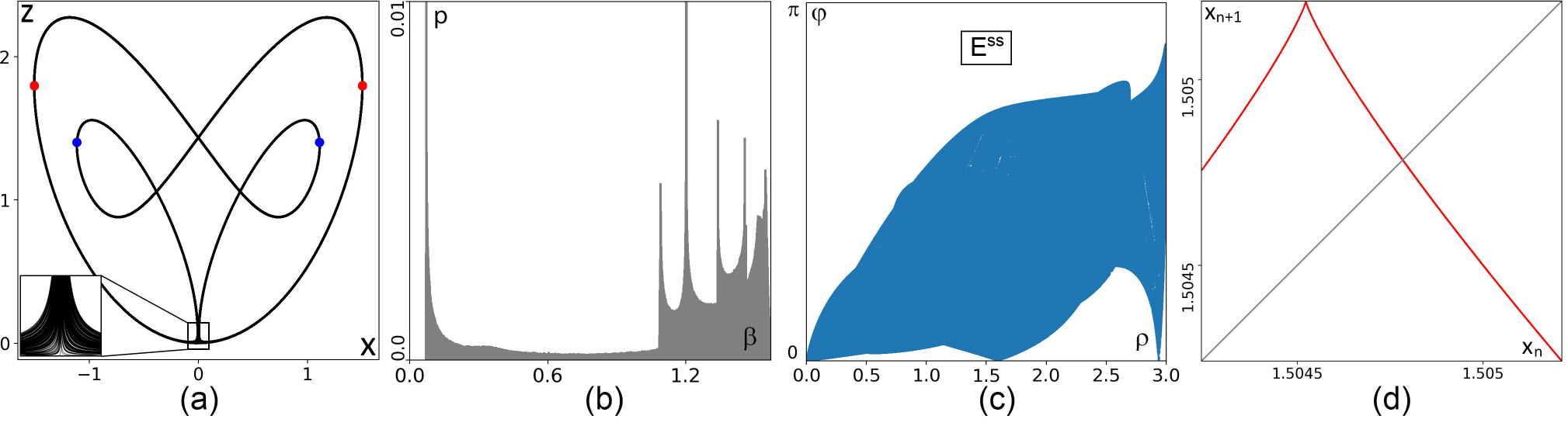}
\caption{\footnotesize Numerical experiments with the Lorenz attractor taken at the point $p_2^1: (\alpha, \lambda) = (0.61, 0.65)$. (a) The projection of the attractor onto the $(x,z)$-plane. The red and blue points correspond to the intersections of the attractor with the cross-section $y=0$ satisfying $\dot z > 0$ (the orbit intersects the cross-section from below). (b) Histogram of the angle between $E^{ss}$ and $E^{cu}$. (c) Continuity diagram for $E^{ss}$. (d) Graph of the 1D Poincar\'e map, numerically constructed from the red intersection points.}
\label{fig_A2_p21}
\end{figure}

\begin{figure}[h!]
\centering
\includegraphics[width=1.0\linewidth]{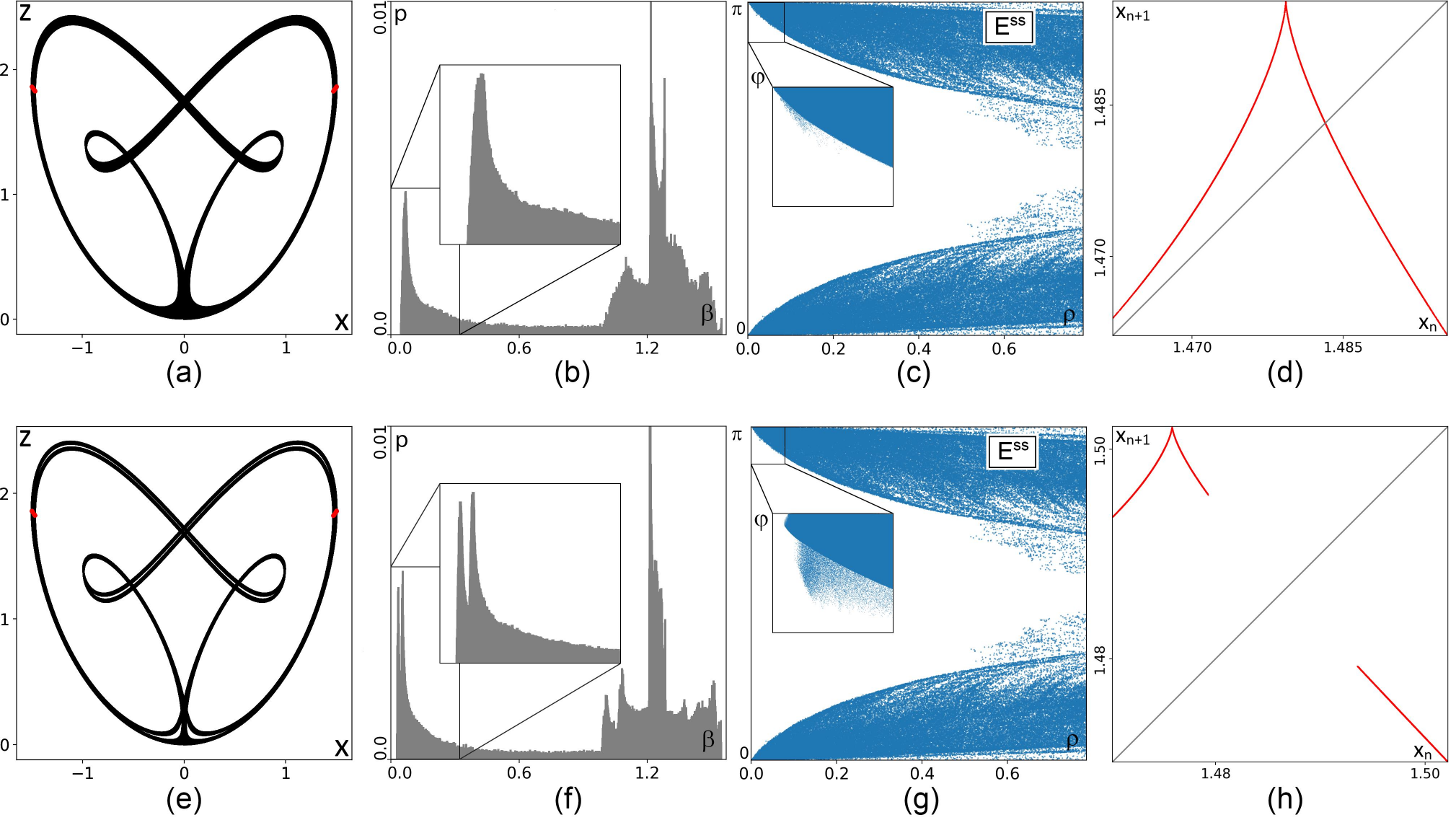}
\caption{\footnotesize Numerical experiments with the Lorenz attractor taken at the point $p_2^2: (\alpha, \lambda) = (0.5, 0.595)$ -- top row; and at the point $p_2^3: (\alpha, \lambda) = (0.508, 0.589)$ -- bottom row. (a), (e) The projection of the attractor onto the $(x,z)$-plane. (b), (f) Histogram of the angle between $E^{ss}$ and $E^{cu}$. (c), (g) Continuity diagram for $E^{ss}$. (d), (h) Graph of the numerically constructed 1D Poincar\'e map.}
\label{fig_A2_p22}
\end{figure}

Now, we describe the changes of the Lorenz attractor near the point $A_2$. In Figures~\ref{fig_A2_p21}--\ref{fig_A2_p22}, we show the results of a numerical experiment for attractors taken at the points $p_2^1$, $p_2^2$, and $p_2^3$ marked in Fig.~\ref{fig_A2}. At the point $p_2^1$, taken above $l_{A=0}$ and near the curve $l_2$, the attractor repeats the shape of the 2-round homoclinic loop and has a large lacuna, see Fig.~\ref{fig_A2_p21}a. The histogram of angles between the subspaces $E^{ss}$ and $E^{cu}$ shown in Fig.~\ref{fig_A2_p21}b confirms that the minimal angle is clearly bounded away from zero, i.e., the attractor is indeed the Lorenz one. The next test establishes the orientability of the attractor. We build the continuity diagram of $E^{ss}$ for a sufficiently long representative orbit of the attractor, as described in Section~\ref{numerical_verification}. The results of this test are shown in Fig.~\ref{fig_A2_p21}c. Since the cloud of points $(\rho_{ij}, \varphi_{ij})$ touches the line $\varphi = 0$ only at one point $(0,0)$, there is no flip of the Lyapunov covariant vector $V_3$. Thus, we claim that the Lorenz attractor under consideration is orientable. 

Additionally, we compute the 1D Poincar\'e map for the considered attractor. For this, the procedure described in Section~\ref{numerical_verification} is slightly modified. We take $\dot x = y = 0$ as a cross-section and consider only the intersections of the separatrix $\Gamma_+$ with this cross-section in the direction of growth of the variable $z$. The red and blue points in Fig.~\ref{fig_A2_p21}a are the traces of the attractor on this cross-section corresponding to the outer and inner turns of $\Gamma_+$ around the equilibria $O_\pm$. Since the attractor has a lacuna and repeats the shape of the 2-round loop, it is convenient to consider every second iteration of the 2D Poincar\'e map for building a representative 1D map. In Fig.~\ref{fig_A2_p21}d, we plot the resulting 1D map based on the outer cloud of points. Note that this map is well fitted with the model map~\eqref{eq_1DMapSigmaTrunc}, cf. Fig.~\ref{fig_A2_p21}d and in Fig.~\ref{fig_ShilCrit_BD}g.

The points $p_2^2$ and $p_2^3$ are taken in the region of the existence of the Lorenz attractor below the curve $l_{A=0}$. The results of the numerical experiments for these points are shown in the top and bottom rows of Fig.~\ref{fig_A2_p22}, respectively. The first column contains the phase portraits of the corresponding attractors. In the second column we present histograms of the angle between the subspaces $E^{ss}$ and $E^{cu}$. We show the continuity diagrams for $E^{ss}$ in the third column and the numerically found 1D Poincar\'e maps in the fourth column. 

The attractor at the point $p_2^2$ is visually similar to the attractor at the point $p_2^1$ (cf. Fig.~\ref{fig_A2_p21}a and Fig.~\ref{fig_A2_p22}a). It is also pseudohyperbolic (the Lorenz attractor), see Fig.~\ref{fig_A2_p22}b. However, this attractor is non-orientable, since the cloud of points in the $E^{ss}$-continuity diagram touches the $\varphi$-axis at two points, $(0,0)$ and $(0,\pi)$, see Fig.~\ref{fig_A2_p22}c. 

A heteroclinic connection occurs on the curve $l_2^{lac}$. The unstable separatrices $\Gamma_\pm$ lie on the stable manifold of the symmetric 4-round saddle periodic orbit. As a result of the heteroclinic bifurcation, the periodic orbit does not belong to the attractor below $l_2^{lac}$ and a new lacuna emerges inside the attractor, see Fig.~\ref{fig_A2_p22}e. The attractor remains the Lorenz attractor, since $\beta > 0$, see Fig.~\ref{fig_A2_p22}f. However, it becomes orientable, since the upper part of the cloud in Fig.~\ref{fig_A2_p22}g is separated from the $\varphi$-axis. The emergence of the lacuna is clearly visible on the numerically constructed 1D Poincar\'e map, cf. Fig.~\ref{fig_A2_p22}d and Fig.~\ref{fig_A2_p22}h.

\begin{figure}[h!]
\centering
\includegraphics[width=1.0\linewidth]{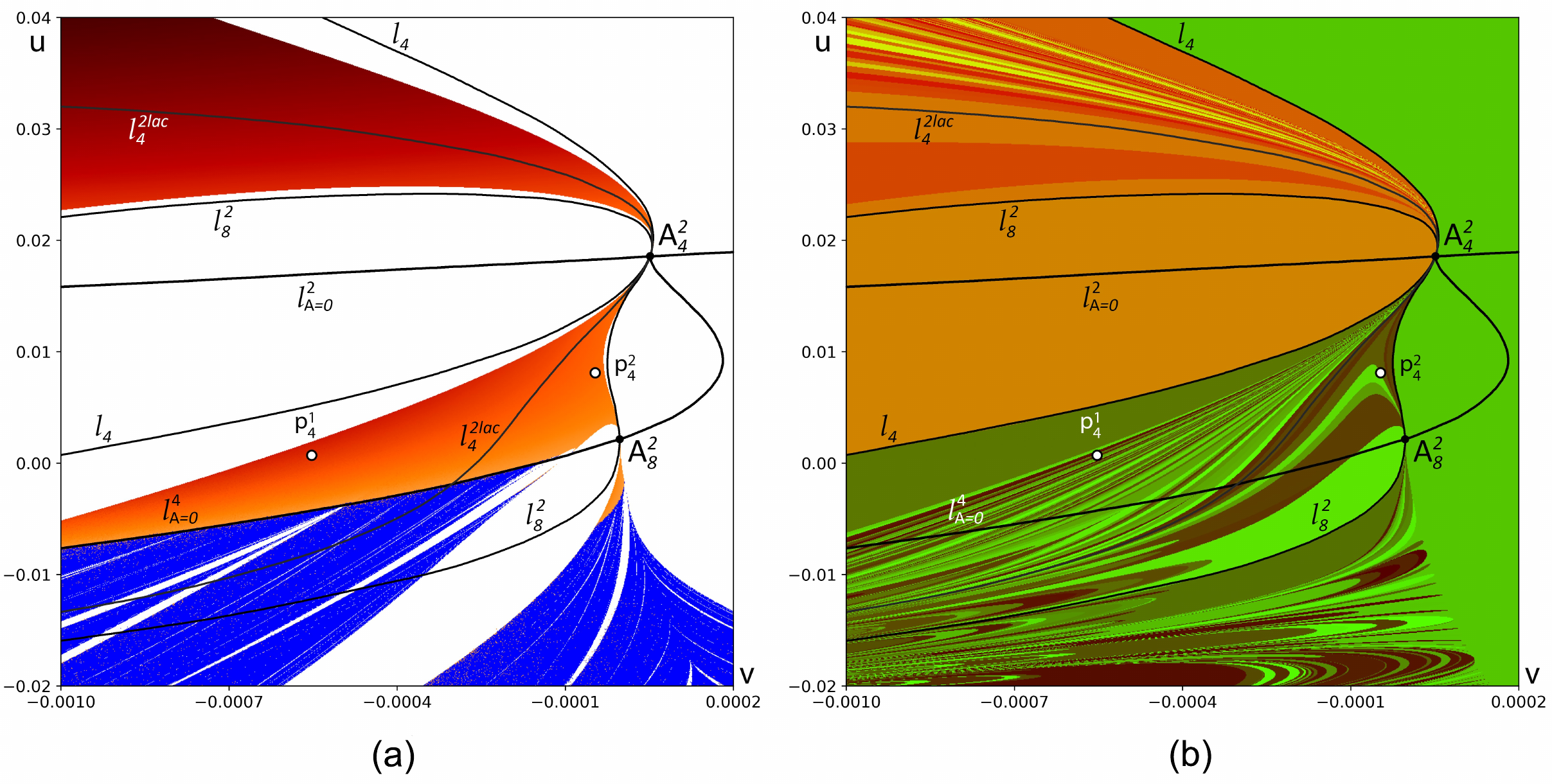}
\caption{\footnotesize Enlarged fragments of (a) the Lyapunov chart and (b) the kneading chart near the point $A_4^2$ in the new parameters $(u,v)$, which are related to the parameters $(\alpha, \lambda)$ by relation \eqref{eq_replace_param}. The curve $l^4_{A=0}$ corresponds to the violation of pseudohyperbolicity condition (P1). This curve starts from the point $A_4^2$ and passes through the point $A_8^2$, where the conditions of the Shilnikov criterion II are satisfied for an 8-round homoclinic loop existing on the curve $l_8^2$. On the curve $l_4^{2lac}$, the unstable separatrices $\Gamma_\pm$ lie on an 8-round symmetric periodic saddle orbit.}
\label{fig_A4}
\end{figure}

\subsection{Cascade of inclination flip points below $l_{A=0}$}\label{sec_cascade_SPII}

The region with the Lorenz attractor shown in Fig.~\ref{fig_A2_p22} is bounded from below by the curve $l^2_{A=0}$, where the second tangency between the subspaces $E^{ss}$ and $E^{cu}$ occurs, see Fig.~\ref{fig_A2}. This curve starts at the point $A_2$ and intersects the homoclinic bifurcation curve $l_4$ at the point $A_4^2$. The point $A_4^2$ corresponds to the inclination flip bifurcation for the pair of 4-round homoclinic loops. The conditions of the Shilnikov criteria II are satisfied at $A_4^2$. 
Therefore, near this point, we observe a similar bifurcation pattern as near $A_2$. In Figure~\ref{fig_A4}, we show an enlarged fragment of Fig.~\ref{fig_A2} near the point $A_4^2$. For panel (a), the threshold of the minimum angle $\beta^*$ is chosen to be equal to $0.00024$. To obtain an informative and noiseless kneading chart in panel (b), we set $K_N = 40$ and skip the first 3 symbols, since this piece of the sequence is always the same. For better clarity, we plot the corresponding chart in the new parameters $u$ and $v$, which are connected with the parameters $\alpha$ and $\lambda$ by the following relation
\begin{equation}
\begin{cases}
    \lambda = 0.48291u - 0.87567v + 0.5704,
    \\
    \alpha = 0.87567u + 0.48291v + 0.47746.
\end{cases}
\label{eq_replace_param}
\end{equation}

\begin{figure}[h]
\centering
\includegraphics[width=1\linewidth]{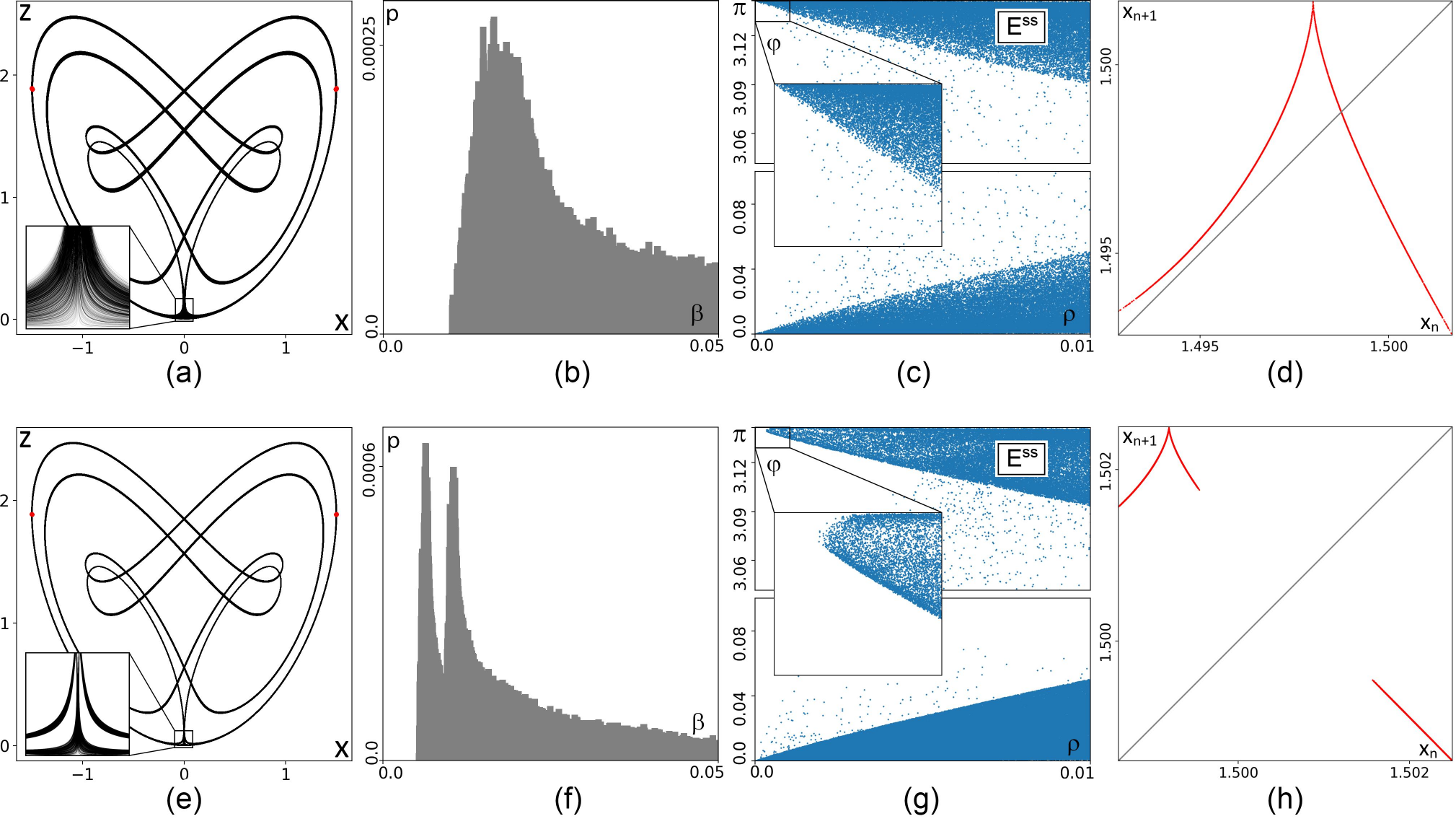}
\caption{\footnotesize Numerical experiments with the Lorenz attractor taken at the point $p_4^1: (u, v) = (-0.0002, -0.0006)$ -- top row; and at the point $p_4^2: (u, v) = (0.01, -0.00005)$ -- bottom row. (a), (e) The projections of the attractor onto the $(x,z)$-plane; (b), (f) Histogram of the angle between $E^{ss}$ and $E^{cu}$; (c), (g) Continuity diagram for $E^{ss}$. (d), (h) Graph of the numerically constructed 1D Poincar\'e map.}
\label{fig_A4_p4}
\end{figure}

Similar to the point $A_2$, two regions with the Lorenz attractor originate from the point $A_4^2$. These regions are enclosed between the curves $l_4$ and $l_8^2$ corresponding to 4- and 8-round homoclinic loops. Nevertheless, the exact boundaries of these regions are formed by heteroclinic connection curves as for the Lorenz attractor born from $A_2$. In Fig.~\ref{fig_A2_p22}e--\ref{fig_A2_p22}h, we show that the attractor at the point $p_2^3$, which is located above the curve $l^2_{A=0}$, is orientable and has two lacuna. In Fig.~\ref{fig_A4_p4}, we provide the results of the numerical experiments for the Lorenz attractor taken at the points $p_4^1$, $p_4^2$, located below the curve $l^2_{A=0}$.

The attractor at the point $p_4^1$ (top row in Fig.~\ref{fig_A4_p4}) also has two lacunae based on 2- and 4-round periodic orbits. This attractor is the Lorenz attractor, since its Lyapunov exponents are equal to $\Lambda_1 = 0.0196, \Lambda_2 = 0, \Lambda_3 = -1.068$ and the minimal angle $\beta_m = 0.0066$ is bounded away from $0$, see Fig.~\ref{fig_A4_p4}b. The continuity diagram of $E^{ss}$ shown in Fig.~\ref{fig_A4_p4}c illustrates that the Lorenz attractor under consideration is non-orientable. In Fig.~\ref{fig_A4_p4}d, we plot the 1D Poincar\'e map which was constructed from the red points lying on the outer part of the attractor. This map looks like the Lorenz map \eqref{eq_1DMapSigmaTrunc}, which additionally confirms that the attractor is the Lorenz one. 

On the curve $l_4^{2lac}$, the unstable separatrices $\Gamma_\pm$ lie on the stable manifold of an 8-round saddle periodic orbit. Passing through this curve, the attractor obtains an additional lacuna based on this periodic orbit. The results of the numerical experiments at the point $p_4^2$, lying to the right of this curve, are given in the bottom row of Fig.~\ref{fig_A4_p4}. The emergence of an additional (third) lacuna is demonstrated by the phase portrait in Fig.~\ref{fig_A4_p4}e and the 1D Poincar\'e map in Fig.~\ref{fig_A4_p4}h. Fig.~\ref{fig_A4_p4}f shows that the attractor remains the Lorenz one. Fig.~\ref{fig_A4_p4}g demonstrates that the attractor becomes orientable, since the cloud of points in the subspace $E^{ss}$ becomes separated from the point $(0,\pi)$.

\begin{figure}[h!]
\centering
\includegraphics[width=1\linewidth]{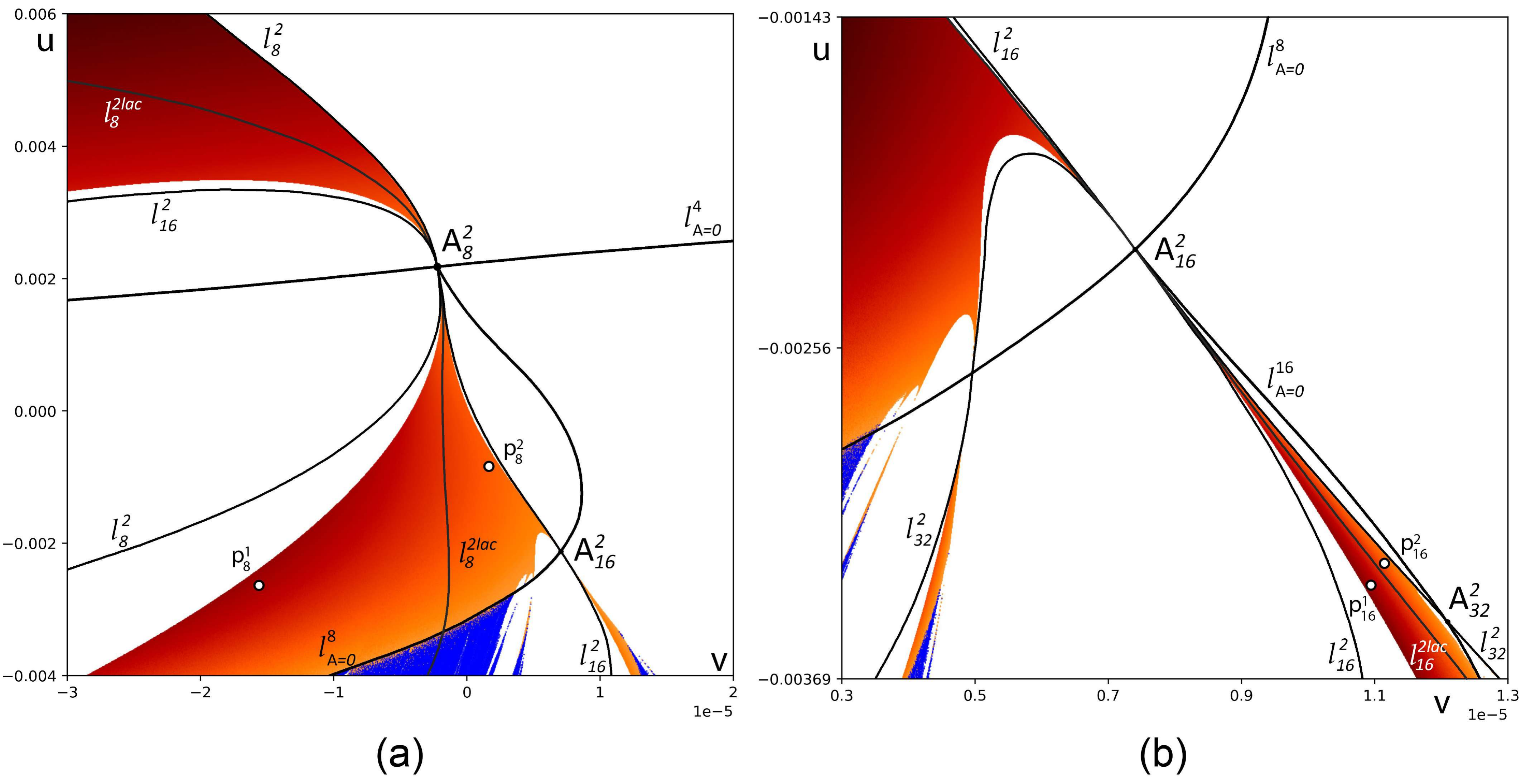}
\caption{\footnotesize Enlarged fragments of the Lyapunov diagram near the Shilnikov points (a) $A_8^2$ and (b) $A_{16}^2$. The curve $l^8_{A=0}$ corresponds to the violation of pseudohyperbolicity condition (P1). This curve starts from the point $A_8^2$ and passes through the point $A_{16}^2$. The 16- and 32-round homoclinic loops occur on the curves $l_{16}^2$ and $l_{32}^2$; on the curves $l_8^{2lac}$ and $l_{16}^{2lac}$ the unstable separatrices $\Gamma_\pm$ lie on the 16- and 32-round symmetric periodic saddle orbit.}
\label{fig_A8_A16}
\end{figure}

In Figure~\ref{fig_A8_A16} we show Lyapunov diagrams near the next inclination flip bifurcation points $A_8^2$, panel (a), and $A_{16}^2$, panel (b). These points are obtained in a similar way as the point $A_8^2$. Namely, the region with the Lorenz attractor shown in Fig.~\ref{fig_A4} is bounded from below by the curve $l^4_{A=0}$ of the tangency between $E^{ss}$ and $E^{cu}$. This curve starts at the point $A_4^2$ and intersects the homoclinic bifurcation curve $l_8^2$ of the 8-round homoclinic loops. The intersection point $A_8^2$ is an inclination flip bifurcation point, where the conditions of the Shilnikov criterion II are satisfied. Similarly, the region with the Lorenz attractor shown in Fig.~\ref{fig_A8_A16}a is bounded from below by the next curve $l^8_{A=0}$ of the tangency between $E^{ss}$ and $E^{cu}$. To plot this curve, we set the threshold of the minimal angle $\beta^*=0.00005$. The curve starts at the point $A_8^2$ and intersects the homoclinic bifurcation curve $l_{16}^2$ of the 16-round homoclinic loops at the point $A_{16}^2$. In Fig.~\ref{fig_A8_A16}b, we also plot the next inclination bifurcation point $A_{32}^2$ of our cascade lying on the curves $l^{16}_{A=0}$ and $l_{32}^2$. To plot the curve $l^{16}_{A=0}$, we set the threshold of the minimal angle $\beta^*=0.00002$.

\begin{figure}[h!]
\centering
\includegraphics[width=1\linewidth]{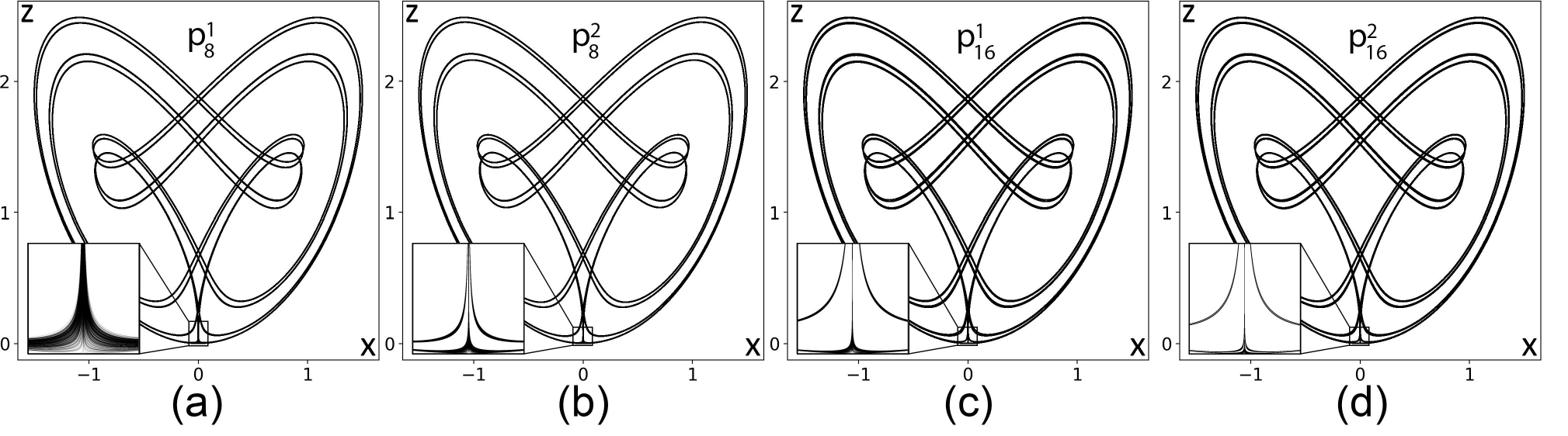}
\vspace{-1cm}
\caption{\footnotesize The phase portrait for the Lorenz attractor taken at the points (a) $p_8^1: (u, v) = (-0.0026, -0.000016)$; (b) $p_8^2: (u, v) = (-0.0007, 0.0000014)$; (c) $p_{16}^1: (u, v) = (-0.0034, 0.000011)$; (d) $p_{16}^2: (u, v) = (-0.003, 0.00001)$.}
\label{fig_pointsA8A16}
\end{figure}

Figure~\ref{fig_pointsA8A16} shows the phase portraits of the attractors taken at the points $p_8^1$ and $p_8^2$ (between $l^4_{A=0}$ and $l^8_{A=0}$) and at the points $p_{16}^1$ and $p_{16}^2$ (between $l^8_{A=0}$ and $l^{16}_{A=0}$). In this figure, we can trace the emergence of additional lacunae associated with the 16- and 32-round periodic orbits. We conjectured that the described cascade of the inclination flip bifurcations is infinite in the Shimizu-Morioka system and leads to the Lorenz attractor with an arbitrarily large number of lacunae. 




\begin{figure}[h!]
\center{\includegraphics[width=0.7\linewidth]{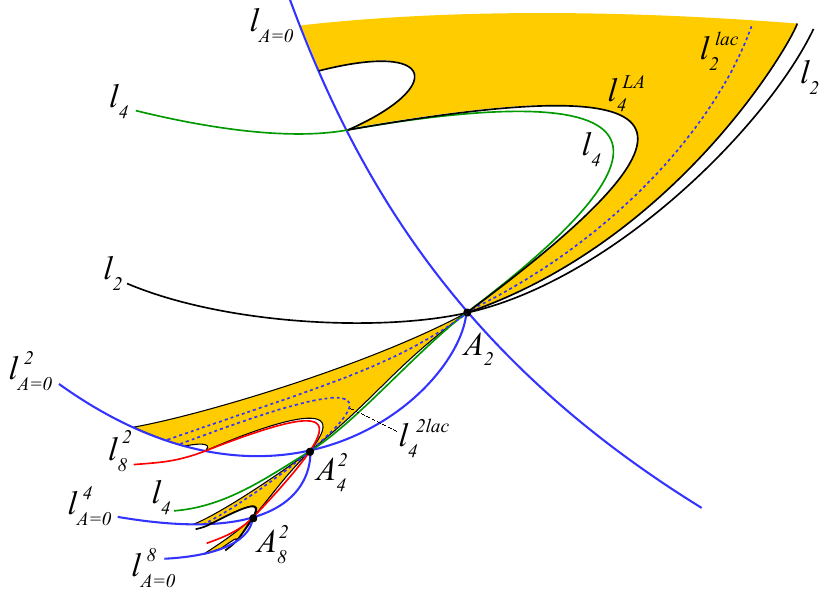}}
\caption{\footnotesize Sketch of the bifurcation diagram illustrating the cascade of the inclination flip bifurcation points. The curves $l^i_{A=0}, i=2, 4, 8, \dots$ correspond to the tangency between $E^{ss}$ and $E^{cu}$. The curves $l^2_i, i=8, 16, \dots$ correspond to $i$-round homoclinic loops. The points $A_2$, and $A^2_i, i = 4, 8, \dots$ are inclination flip bifurcation points, where the conditions of the Shilnikov criterion II are satisfied.}
\label{fig_SM_SketchA0_2}
\end{figure}

To conclude our considerations, we present a schematic illustration of the bifurcation diagram for the cascade described above, see Fig.~\ref{fig_SM_SketchA0_2}. The curves $l_{A=0}$ and $l^i_{A=0},\, i = 2, 4, \dots$ are curves of the tangency between $E^{ss}$ and $E^{cu}$. The curves $l_2$, $l_4$, and $l_i^2,\, i=8, 16, \dots$ are homoclinic bifurcation curves. The points $A_2$, and $A^2_i, i = 4, 8, \dots$ are inclination flip bifurcation points, where the conditions of the Shilnikov criterion II are satisfied. Each curve $l^i_{A=0}, \, i = 2, 4, \dots$ begins at the point $A^2_i$ and intersects the curve $l^2_{2i}$ at the next point $A^2_{2i}$. The Lorenz attractor exists in the yellow region. This region is bounded by heteroclinic bifurcation curves where $\Gamma_{\pm}$ lies on the stable manifolds of saddle periodic orbits. We note that along each curve $l^i_{A=0}$, the region of the existence of the Lorenz attractor should have the same structure as near the first curve $l_{A=0}$, which was described in Section~\ref{sec_cascade_SPI}.

We also note that for an appropriate one-parameter pathway crossing the diagram in Fig.~\ref{fig_SM_SketchA0_2}, we can observe the implementation of the phenomenological scenario described in Section~\ref{sec_cascade}. The intersections of this pathway with the curves $l_2$, $l_4$, and $l_i^2,\, i=8, 16, \dots$ correspond (on a properly chosen cross-section) to $M_{\pm}$ lying on the line $\Pi_{0}$, as show in Figs.~\ref{fig_GM1a}a, \ref{fig_GM2}a, and \ref{fig_GM2}d. The intersections with the yellow region of the existence of the Lorenz attractor correspond to $M_{\pm}$ lying on the lines $\pi^i_{\pm}$, as shown in Figs.~\ref{fig_GM1a}c, \ref{fig_GM2}c, and \ref{fig_GM2}f. New hooks in the wedges $T(\Pi_{\pm})$ occur when the one-parameter way intersects the curves $l^i_{A=0}$.

\section*{Acknowledgements}
This work was supported by the Basic Research Program at HSE University.

\printbibliography

\end{document}